\documentclass[10pt,leqno,draft,twoside]{article}

\usepackage{amsmath,amsthm,amssymb}
\numberwithin{equation}{section}

\makeatletter
\newif\ify@autoscale \y@autoscaletrue \def\Yautoscale#1{\ifnum #1=0
  \y@autoscalefalse\else\y@autoscaletrue\fi}
\newdimen\y@b@xdim
\newdimen\y@boxdim \y@boxdim=13pt
\def\y@style{\displaystyle} 
\def\Yboxdim#1{\y@autoscalefalse\y@boxdim=#1}
\newdimen\y@linethick    \y@linethick=.3pt
\def\Ylinethick#1{\y@linethick=#1}
\newskip\y@interspace \y@interspace=0ex plus 0.3ex
\def\Yinterspace#1{\y@interspace=#1}
\newif\ify@vcenter   \y@vcentertrue	
\def\Yvcentermath#1{\ifnum #1=0 \y@vcenterfalse\else\y@vcentertrue\fi}
\newif\ify@stdtext   \y@stdtextfalse
\def\Ystdtext#1{\ifnum #1=0 \y@stdtextfalse\else\y@stdtexttrue\fi}
\def\y@vr{\vrule height0.8\y@b@xdim width\y@linethick depth 0.2\y@b@xdim}
\def\y@emptybox{\y@vr\hbox to \y@b@xdim{\hfil}}
\def\y@abcbox#1{\y@vr\hbox to \y@b@xdim{\hfil#1\hfil}}
\def\y@mathabcbox#1{\y@vr\hbox to \y@b@xdim{\hfil${\y@style#1}$\hfil}} 
\def\y@setdim{%
  \ify@autoscale%
   \setbox1=\hbox{A}\y@b@xdim=1.6\ht1 \setbox1=\hbox{}\box1%
  \else\y@b@xdim=\y@boxdim \advance\y@b@xdim by -2\y@linethick
  \fi}
\newcount\y@counter
\newif\ify@islastarg
\def\y@lastargtest#1,#2 {\if\space #2 \y@islastargtrue
  \else\y@islastargfalse\fi}
\def\y@emptyboxes#1{\y@counter=#1\loop\ifnum\y@counter>0
  \advance\y@counter by -1 \y@emptybox\repeat}
\def\y@nelineemptyboxes#1{%
  \vbox{%
    \hrule height\y@linethick%
    \hbox{\y@emptyboxes{#1}\y@vr}
    \hrule height\y@linethick}\vspace{-\y@linethick}}
\def\yng(#1){%
  \y@setdim%
  \hspace{\y@interspace}%
  \ifmmode\ify@vcenter\vcenter\fi\fi{%
  \y@lastargtest#1,
  \vbox{\offinterlineskip
    \ify@islastarg
     \y@nelineemptyboxes{#1}
    \else
     \y@ungempty(#1)
    \fi}}\hspace{\y@interspace}}
\def\y@ungempty(#1,#2){%
  \y@nelineemptyboxes{#1}
  \y@lastargtest#2,
  \ify@islastarg
   \y@nelineemptyboxes{#2}
  \else
   \y@ungempty(#2)
  \fi}
\def\y@nelettertest#1#2. {\if\space #2 \y@islastargtrue
  \else\y@islastargfalse\fi}
\def\y@abcboxes#1#2.{%
  \ify@stdtext\y@abcbox#1\else\y@mathabcbox#1\fi%
  \y@nelettertest #2.
  \ify@islastarg\unskip%
   \ify@stdtext\y@abcbox{#2}\else\y@mathabcbox{#2}\fi%
  \else\y@abcboxes#2.\fi}
\def\y@m@veright@ifskew#1{}
\def\y@nelineabcboxes#1{%
  \y@nelettertest #1.
  \ify@islastarg
   \y@m@veright@ifskew{#1}
    \vbox{
      \hrule height\y@linethick%
      \hbox{\ify@stdtext\y@abcbox#1\else\y@mathabcbox#1\fi\y@vr}
      \hrule height\y@linethick}\vspace{-\y@linethick}
  \else
   \y@m@veright@ifskew{#1}
    \vbox{
      \hrule height\y@linethick%
      \hbox{\y@abcboxes #1.\y@vr}%
      \hrule height\y@linethick}\vspace{-\y@linethick}
  \fi}
\def\young(#1){%
  \y@setdim%
  \hspace{\y@interspace}%
  \y@lastargtest#1,
  \ifmmode\ify@vcenter\vcenter\fi\fi{%
  \vbox{\offinterlineskip
    \ify@islastarg\y@nelineabcboxes{#1}%
    \else\y@ungabc(#1)%
    \fi}}\hspace{\y@interspace}}
\def\y@ungabc(#1,#2){%
  \y@nelineabcboxes{#1}%
  \y@lastargtest#2,
  \ify@islastarg\y@nelineabcboxes{#2}%
  \else\y@ungabc(#2)%
  \fi}

\def\minitab(#1){%
\y@autoscalefalse%
\def\y@style{\scriptscriptstyle}%
\y@boxdim=6pt%
\young(#1)%
\y@autoscaletrue%
\def\y@style{\displaystyle}%
\y@boxdim=13pt}
\def\tinytab(#1){%
\y@autoscalefalse%
\def\y@style{\scriptscriptstyle}%
\y@boxdim=5pt%
\young(#1)%
\y@autoscaletrue%
\def\y@style{\displaystyle}%
\y@boxdim=13pt}

\def\miniyng(#1){%
\y@autoscalefalse%
\y@boxdim=6pt%
\yng(#1)%
\y@autoscaletrue%
\y@boxdim=13pt}
\def\tinyyng(#1){%
\y@autoscalefalse%
\y@boxdim=5pt%
\yng(#1)%
\y@autoscaletrue%
\y@boxdim=13pt}
\def\nanoyng(#1){%
\y@autoscalefalse%
\y@boxdim=3pt%
\yng(#1)%
\y@autoscaletrue%
\y@boxdim=13pt}

\makeatother

\newcommand{\by}[2]{#1, #2}  
\newcommand{\journal}[1]{{\itshape #1}}
\newcommand{\volume}[1]{{\bfseries #1}}
\newcommand{\doi}[1]{\texttt{doi:#1}}

\newtheorem{thm}{Theorem}[section]
\newtheorem{prop}[thm]{Proposition}
\newtheorem{lem}[thm]{Lemma}
\newtheorem{cor}[thm]{Corollary}
\newtheorem{conj}[thm]{Conjecture}

\theoremstyle{definition}

\newtheorem{dfn}[thm]{Definition}
\newtheorem{example}[thm]{Example}

\theoremstyle{remark}
\newtheorem{remark}[thm]{Remark}

\newcommand{\eoe}{\hfill{\rule{0.7em}{0.7em}}}

\newenvironment{ex}{\begin{example}}
{\eoe\end{example}}
\newenvironment{rem}{\begin{remark}}
{\eoe\end{remark}}

\newcommand{\deq}{:=}

\newcommand{\card}[1]{\left|#1\right|}
\newcommand{\abs}[1]{\left|#1\right|}
\newcommand{\sym}[1]{\mathfrak{S}_{#1}}

\newcommand{\N}{\mathbb{N}}
\newcommand{\Z}{\mathbb{Z}}

\newcommand{\C}{\mathbb{C}}

\newcommand{\multsum}[1]{\sum_{\substack{#1}}}
\newcommand{\multprod}[1]{\prod_{\substack{#1}}}

\newcommand{\then}{\,\Longrightarrow\,}

\newcommand{\set}[3][0]
{\ifcase#1%
\left\{#2\,;\,#3\right\}\or%
\bigl\{#2\,\big|\,#3\bigr\}\or%
\Bigl\{#2\,\Big|\,#3\Bigr\}\or%
\biggl\{#2\,\bigg|\,#3\biggr\}\or%
\Biggl\{#2\ \Bigg|\ #3\Biggr\}\else%
\{#2;#3\}\fi}

\newcommand{\genkakko}[4]
{\ifcase#1%
\left#3#2\right#4\or%
\bigl#3#2\bigr#4\or%
\Bigl#3#2\Bigr#4\or%
\biggl#3#2\biggr#4\or%
\Biggl#3#2\Biggr#4\else%
#2\fi}

\newcommand{\kakko}[2][0]{\genkakko{#1}{#2}{(}{)}}
\newcommand{\ckakko}[2][0]{\genkakko{#1}{#2}{\{}{\}}}

\newcommand{\inprod}[3][0]{\genkakko{#1}{#2,\,#3}{\langle}{\rangle}}

\newcommand{\vkdet}[2][k]{\left|#2\right|_{#1}}
\newcommand{\pmat}[1]{P(#1)}

\let \det = \relax
\DeclareMathOperator{\det}{det}
\DeclareMathOperator{\Imm}{Imm}
\DeclareMathOperator{\Mat}{Mat}
\DeclareMathOperator{\sgn}{sgn}
\DeclareMathOperator{\STab}{STab}
\DeclareMathOperator{\SSTab}{SSTab}
\DeclareMathOperator{\ML}{ML}
\DeclareMathOperator{\AL}{AL}
\DeclareMathOperator{\tdet}{det}
\DeclareMathOperator{\diag}{diag}

\DeclareMathOperator{\ind}{ind}
\DeclareMathOperator{\res}{res}

\newcommand{\liegl}{\mathfrak{gl}}

\newcommand{\A}{\mathcal{P}} 
\newcommand{\U}{\mathcal{U}}

\newcommand{\va}{\boldsymbol{a}}
\newcommand{\vb}{\boldsymbol{b}}
\newcommand{\ve}{\boldsymbol{e}}

\newcommand{\vj}{\boldsymbol{j}}

\newcommand{\1}{\boldsymbol{1}}
\newcommand{\vdelta}{\boldsymbol{\delta}}
\newcommand{\Uq}{\mathcal{U}_q(\mathfrak{gl}_n)}

\newcommand{\trivrpn}[1]{1_{#1}}

\newcommand{\tp}[1]
{{}^{\raise0.35ex\hbox{$\scriptstyle t$}}\kern-0.2em #1}

\newcommand{\xs}{x_{\star}}

\newcommand{\kple}[2][k]{{#2}^{[#1]}}
\newcommand{\kplerow}[2][k]{{#2}_{[#1]}}
\newcommand{\kplerc}[2][k]{{#2}^{[#1]}_{[#1]}}
\newcommand{\adet}[1][\alpha]{\det^{(#1)}\!}
\newcommand{\kdet}[1][k]{\det_{#1}}
\newcommand{\wrdet}[1][k]{\operatorname{wrdet}_{#1}}
\newcommand{\len}[1]{\ell(#1)}
\newcommand{\ksgn}[1][n,k]{\operatorname{sgn}_{#1}}

\newcommand{\wrprod}[2]{#1\wr\sym#2}
\newcommand{\wsym}[2]{\wrprod{\sym#1}{#2}}

\newcommand{\GLmod}[2]{{\mathcal M}_{#1}^{#2}}
\newcommand{\Smod}[2]{{\mathcal J}_{#1}^{#2}}

\newcommand{\rnk}{\mathfrak{R}}
\newcommand{\floor}[1]{\left\lfloor#1\right\rfloor}
\newcommand{\symp}{\mathcal{S}}

\DeclareMathOperator{\wt}{wt}

\newcommand{\pfkakko}[1]{{\upshape(}#1\/{\upshape)}}
\newcommand{\pfcite}[1]{{\upshape\cite{#1}}}
\newcommand{\pfref}[1]{{\upshape\ref{#1}}}

\newcommand{\snk}{\mathfrak{R}^{\times}}

\newenvironment{question}{\medskip\itshape}{\medskip}
\newenvironment{keywords}{\smallskip\noindent{\bfseries Keywords:}}{}
\newenvironment{MSC}{\smallskip\noindent{\bfseries 2000 Mathematical Subject Classification:}}{}
\newenvironment{dedication}{\begin{center}\itshape}{\end{center}}

\begin{document}

\title{\bfseries
Invariant theory for singular $\alpha$-determinants}
\author{Kazufumi Kimoto%
\thanks{Partially supported by Grant-in-Aid for Young Scientists (B) No.16740021.}
\and
Masato Wakayama%
\thanks{Partially supported by Grant-in-Aid for Scientific Research (B) No. 15340012.}}
\date{\today}
\pagestyle{myheadings}
\markboth{Invariant theory for singular $\alpha$-determinants}
{K. Kimoto and M. Wakayama.}

\maketitle

\begin{dedication}
Dedicated to Roger Howe on the occasion of his 60th birthday
\end{dedication}

\begin{abstract}
From the irreducible decompositions' point of view,
the structure of the cyclic $GL_n(\C)$-module
generated by the $\alpha$-determinant degenerates when
$\alpha=\pm \frac1k \;(1\leq k\leq n-1)$ (see \cite{MW2005}).
In this paper, we show that $-\frac1k$-determinant
shares similar properties which the ordinary determinant
possesses. From this fact, one can define a new (relative)
invariant called a \emph{wreath determinant}.
Using $(GL_m,\, GL_n)$-duality
in the sense of Howe, 
we obtain an
expression of a wreath determinant by a certain
linear combination of the corresponding ordinary minor determinants
labeled by suitable rectangular shape tableaux.
Also we study a wreath determinant analogue of the
Vandermonde determinant, and then, investigate
symmetric functions such as Schur functions
in the framework of wreath determinants.
Moreover, we examine coefficients which we call $(n,k)$-sign
appeared at the linear expression of the wreath determinant
in relation with a zonal spherical function of a Young subgroup
of the symmetric group $\sym {nk}$.


\begin{MSC}
17B10, 
15A72, 
05E10. 
\end{MSC}
\end{abstract}

\tableofcontents

\section{Away from the multiplication law}

There is a notion called the \emph{$\alpha$-determinant}
for a square matrix in probability theory.
It was first introduced in \cite{V}
and actually appeared as coefficients of the Taylor expansion
of $\det(I-\alpha A)^{-1/\alpha}$.
This expansion has applications,
in particular, to multivariate binomial and negative binomial distributions.
Moreover, recently in \cite{ST2003},
the $\alpha$-determinant is use to define a random point process
through a study of the Fredholm determinants of certain integral operators.

The $\alpha$-determinant $\adet(X)$ for a matrix $X$
(see \eqref{eq:def_of_alpha_det} for the definition)
does not have the multiplication property
which the ordinary determinant $\det(X)$ possesses.
It is, however, interesting from a viewpoint of invariant theory
because the $\alpha$-determinant is regarded as an interpolation
of the determinant ($\alpha=-1$) and permanent ($\alpha=1$)
--- recall that each of them generates
an irreducible representation of the general linear group $GL_n(\C)$;
as representations of the special linear group $SL_n(\C)$,
the former defines the trivial representation
and the latter generates the representation
on the space of symmetric $n$-tensors
of (the natural representation on) $\C^n$.
These facts raise naturally the following question:

\begin{question}
``Where had the multiplication law gone
when $\alpha$ moved away from $-1$?''
\end{question}

The multiplication law of the determinant is equivalent to the fact that
$GL_n(\C)\cdot \det(X)\subset \C^\times \det(X)$. Hence, it is
natural to ask the question what the smallest invariant space
containing $GL_n(\C)\cdot \adet(X)$ is.
From this point of view,
Matsumoto and the second author \cite{MW2005} have studied recently
the irreducible decomposition of the cyclic module
$\U(\liegl_n)\cdot\adet(X)$
and showed that
the structure of the module changes drastically
when $\alpha$ is contained in the set
$\{\pm1,\pm\frac12,\dots,\pm\frac1{n-1}\}$.
In fact, one can see that
the irreducible decomposition of the cyclic module
$\U(\liegl_n)\cdot\adet(X)$ degenerates
when $\alpha$ is one of such values. More precisely,  if we denote by
$m^\lambda(\alpha)$ the multiplicity of the irreducible 
highest weight  $\U(\liegl_n)$-module 
corresponding to a partition 
$\lambda$ appeared in the decomposition,
then, for instance, we have $m^\lambda(-\frac1k)=0$ when the first
component of $\lambda$ is greater than $k$ 
(see \eqref{eq:degeneration_of_cyclic_module}).
Therefore, we shall call $\alpha$ \emph{singular} if
$\alpha\in\{\pm1,\pm\frac12,\dots,\pm\frac1{n-1}\}$.
This result indicates that
if $\alpha$ is singular, then
$\adet(X)$ may share some distinguished feature
which explains
why such a drastic change of the module structure happens.
The special emphasis in this paper is laid on the study 
of the case $\alpha=-\frac1k\; (k\in \Z_{>0})$. 
Actually, we first show that $\det^{(-\frac1k)}(X)$ has a certain
alternating property which is considered as a generalization
of the alternating property of the ordinary determinant
(as well as its multilinearity) in Section \ref{sec:lemmas}.
We also show that such an alternating property
characterizes the $-\frac1k$-determinants
through the cyclic module $\U(\liegl_n)\cdot\det^{(-\frac1k)}(X)$
by the effective use of the Young symmetrizer
(Section \ref{sec:characterization}).
We note that a quantum analogue of the $\alpha$-determinant
(which we call quantum $\alpha$-determinant)
is introduced and studied in \cite{KW2007},
however, it is much more difficult to describe
the singular values in the quantum case.

Under these studies, one of the main purpose of the present paper
is to construct an invariant,
which we will call a \emph{wreath determinant},
defined by means of a singular $\alpha$-determinant.
In order to obtain this new invariant for a rectangular matrix,
we consider a $kn\times kn$ matrix gotten from multiplexing
a given $kn\times n$ matrix $A$
by tensoring the $1\times k$ matrix $(1,1,\ldots,1)$. 
By using the property of $\alpha$-determinants developed in
Section \ref{sec:lemmas} for $\alpha=-\frac1k$, we show that the wreath 
determinant is a relative invariant for the action of the 
wreath product of symmetric groups $\wsym kn$ (see \cite{Mac}) in 
Section \ref{sec:multiple}.
Furthermore, in Section \ref{sec:detexp},
we give an expression of the wreath determinant of
$kn\times n$-matrix $A$ by a linear sum of the
$n$-th minor determinants of $A$ labeled by the corresponding
rectangular shaped tableaux.
In the derivation of this expression,
$(GL_m,\,GL_n)$-duality in the sense of \cite{Howe1}
provides a guiding principle.
We then, beside the expression above, derive another
expression of such a wreath determinant conceptually 
by the Frobenius reciprocity.
As a corollary of the proof, we find that 
the wreath determinant is a relative invariant of 
$(\wsym kn)\times GL_n$.
We also give one remark on the background which explains how to 
get this expression and to understand a structure of the cyclic module 
$\U(\liegl_n)\cdot\adet(X)^\ell$ for a general 
positive integer $\ell$ 
in the framework of $(GL_m,\,GL_n)$-duality.
Note that the latter closely  
relates a problem for calculating a certain plethysm
\cite{Mac, Howe2}.

The Cauchy determinant formula (see, e.g. \cite{Weyl})
\begin{equation*}
\det\kakko{\frac1{x_i+y_j}}_{1\le i,j\le n}
=\frac{\Delta_n(x)\Delta_n(y)}{\prod_{i,j=1}^n(x_i+y_j)}
\end{equation*}
can be considered as one of the most important determinant formula
from the representation theoretic point of view.
In Section \ref{sec:cauchy},
we prove an analogue of
the Cauchy determinant formula for the wreath determinants.
It naturally leads us to study the wreath determinant of a 
Vandermonde type. 
The aforementioned study enables us to deduce
a formula for the Schur functions in terms of the $-\frac1k$-determinants
of the Vandermonde type,
which is regarded as a $-\frac1k$-analogue of
the expression
\begin{equation*}
s_\lambda(x_1,\dots,x_n)
=\frac{\det(x_i^{\lambda_j+n-j})_{1\le i,j\le n}}
{\det(x_i^{n-j})_{1\le i,j\le n}}.
\end{equation*}
The proof is to be done first for the corresponding expressions
for the monomial symmetric functions $m_\lambda(x)$,
and then, it can be completed immediately
by the well-known linear expression of the
Schur function by  $m_\lambda(x)$ using the Kostka numbers
(Section \ref{sec:cauchy}).

We further try to understand the coefficients
which we call \emph{$(n,k)$-sign}
appeared at the aforementioned
linear expression of the wreath determinant
in relation with a zonal spherical function of
a Young subgroup of the symmetric group $\sym {nk}$.
At this point, we shall provide one conjecture about
a positive definiteness of a certain symmetric matrix formed by the
spherical function (see Conjecture \ref{positiveConjecture}).
We don't treat the remaining singular case 
$\alpha=\frac1k\; (k\in \Z_{>0})$. Note that, however, one can deduce 
the fact $m^\lambda(\frac1k)=m^{\lambda'}(-\frac1k)$ from the result in 
\cite{MW2005}, where $\lambda'$ denotes the transposition of 
the partition $\lambda$ as a Young diagram. 

We give an $\alpha$-analogue of the Laplace expansion formula
for $\alpha$-determinants in Appendix \ref{sec:laplace}.

\subsection*{Conventions}

As usual, $\N$ is the set of positive integers and
$\C$ is the complex number field.
For $n\in\N$, we denote by $\sym n$ the symmetric group of degree $n$.
The cycle number of an element $\sigma\in\sym n$ is written by $\nu_n(\sigma)$.
Since the conjugacy classes of $\sym n$ are parametrized by the cycle type,
$\nu_n$ is a class function on $\sym n$.
In particular, we notice that
$\nu_n(\sigma^{-1})=\nu_n(\sigma)$ for any $\sigma\in\sym n$
because $\sigma$ and $\sigma^{-1}$ are always $\sym n$-conjugate.

We denote by $\Mat_{m,n}$ the set of $m\times n$ matrices
whose entries belong to a certain commutative $\C$-algebra,
and we put $\Mat_n=\Mat_{n,n}$.
We also denote by $I_n=(\delta_{ij})_{1\le i,j\le n}$ the identity matrix of size $n$
and $\1_n=(1)_{1\le i,j\le n}$ the all-one matrix of size $n$.
For a permutation $\sigma\in\sym n$,
$\pmat\sigma=(\delta_{i\sigma(j)})_{1\le i,j\le n}$
is the permutation matrix for $\sigma$.

The (complex) general linear group $GL_n(\C)$ is the group
consisting of invertible matrices in $\Mat_n(\C)$.
We exclusively deal with the complex vector spaces
so that we often omit the symbol $\C$ and
simply write $GL_n$ in stead of writing $GL_n(\C)$.

Let us put $[N]\deq\{1,2,\dots,N\}$ for $N\in\N$.
For a given partition (or Young diagram) $\lambda$ of size $N$,
we denote by $\SSTab_N(\lambda)$
the set of all semistandard tableaux
with shape $\lambda$ whose entries are in $[N]$,
and we also denote by $\STab(\lambda)$
the set of all standard tableaux with shape $\lambda$.
For a semistandard tableau $T\in\SSTab_N(\lambda)$,
we associate a sequence $\wt(T)\deq(\mu_1,\mu_2,\dots,\mu_N)$
of nonnegative integers where $\mu_k=\card{\{t_{ij}=k\}}$
is the number of entries in $T$ which is equal to $k$.
We call $\wt(T)$ the \emph{weight} of $T$.
Notice that a semistandard tableau $T\in\SSTab_N(\lambda)$ is standard
if and only if $\wt(T)=(1,1,\dots,1)$.
For a given partition $\lambda,\mu\vdash N$ of the same size $N$,
we denote by $K_{\lambda\mu}$ the number of semistandard tableaux $T$
with shape $\lambda$ such that $\wt(T)=\mu$.
Namely,
\begin{equation*}
K_{\lambda\mu}=\card{\set[1]{T\in\SSTab_N(\lambda)}{\wt(T)=\mu}}.
\end{equation*}
We call $K_{\lambda\mu}$ the \emph{Kostka number}.
We also put $f^\lambda=\card{\STab(\lambda)}=K_{\lambda,(1,\dots,1)}$,
and denote by $\ell(\lambda)$ the depth of the diagram $\lambda$.
See \cite{F, Mac} for detailed information on partitions and tableaux.

The irreducible polynomial representations of $GL_m$
are highest weight modules and the highest weights are identified
with partitions such that $\len\lambda\le m$.
We denote by $\GLmod m\lambda$ the irreducible $GL_m$-module
corresponding to the partition $\lambda$.
The irreducible representations of $\sym n$ are
also parametrized by partitions of $n$.
We denote by $\Smod n\lambda$ the irreducible $\sym n$-module
corresponding to the partition $\lambda\vdash n$.
See \cite{Weyl} (or \cite{F}) for detailed information on representation theory
of $GL_m$ and $\sym n$.

\section{Basic properties of general $\alpha$-determinants}\label{sec:lemmas}

Let $\alpha$ be a complex parameter.
The $\alpha$-determinant $\det^{(\alpha)}A$ of
a square matrix $A=(a_{ij})_{1\le i,j\le n}\in\Mat_n$ is defined by
\begin{equation}\label{eq:def_of_alpha_det}
\adet A\deq\sum_{w\in\sym n}\alpha^{n-\nu_n(w)}a_{w(1)1}\dotsb a_{w(n)n}.
\end{equation}
We note that $\adet{}(\tp{A})=\adet(A)$
because $\nu_n(w^{-1})=\nu_n(w)$ for any $w\in\sym n$.
We also notice that $\adet$ is \emph{multilinear}
with respect to the column (and/or row) vectors.
We mainly deal with the $-\frac1k$-determinants
for $k\in\N$ below,
so it is convenient to put
\begin{equation*}
\kdet A=\vkdet A\deq\adet[-1/k] A.
\end{equation*}
We note that $\kdet[1]=\adet[-1]$ is the ordinary determinant.

The $\alpha$-determinant of the all-one matrix $\1_n$ (i.e. every element equals $1$) is
calculated as
\begin{align}\label{eq:adet_of_1}
\adet\1_n=\sum_{w\in\sym n}\alpha^{n-\nu_n(w)}
=\prod_{1\le i<n}(1+i\alpha).
\end{align}
We note that this is the generating function of
the Stirling numbers of the first kind
(see, e.g. \cite{St}).
The following lemma is
the (shifted) partial sum generalization of the identity above.

\begin{lem}\label{lem:shifted_cycle_sum}
For a subset $I$ of $[n]=\{1,2,\dots,n\}$, put
\begin{equation*}
\sym n(I)\deq
\set[2]{w\in\sym n}{x\notin I \then w(x)=x}.
\end{equation*}
Then, for any $g\in\sym n$,
there exists a nonnegative integer $m(g,I)$ such that
\begin{equation*}
\sum_{w\in\sym n(I)}\alpha^{n-\nu_n(gw)}
=\alpha^{m(g,I)}\prod_{1\le i<k}(1+i\alpha),
\end{equation*}
where $k=\card{I}$.
The integer $m(g,I)$ is given by $n-\nu_n(gw_0)$
where $w_0\in\sym n(I)$ is the unique element
such that $\nu_n(gw_0)\ge\nu_n(gw)$
for any $w\in\sym n(I)$.
\end{lem}

\begin{proof}
Take an element $h\in\sym n$ such that $h\cdot I=[k]$.
We identify $\sym k$ and $\sym n([k])$ naturally.
Since $w\in\sym n(I)$ if and only if $hwh^{-1}\in\sym k$,
it follows that
\begin{equation*}
\sum_{w\in\sym n(I)}\alpha^{n-\nu_n(gw)}
=\sum_{w\in\sym k}\alpha^{n-\nu_n(gh^{-1}wh)}
=\sum_{w\in\sym k}\alpha^{n-\nu_n(gh^{-1}(hw_0h^{-1})wh)}
=\sum_{w\in\sym k}\alpha^{n-\nu_n(g'w)}
\end{equation*}
where $g'=hgw_0h^{-1}$.
By the definition of $w_0$ and $g'$,
it is easy to see that
\begin{equation}\label{eq:nu-ineq}
\nu_n(g')\ge\nu_n(g'w)
\quad(w\in\sym k).
\end{equation}
Assume that
$g'$ contains a cycle of the form $(\vj_2,i_2,\vj_1,i_1)$
($i_1,i_2\in\{1,2,\dots,k\}$, $i_1\ne i_2$ and $\vj_1,\vj_2$ stand for
certain disjoint strings
in $\{1,2,\dots,n\}$ which are possibly empty).
Then it follows that $\nu_n(g'\cdot(i_1,i_2))=\nu_n(g')+1$
because
\begin{equation*}
(\vj_2,i_2,\vj_1,i_1)\cdot(i_2,i_1)
=(\vj_2,i_2)\cdot(\vj_1,i_1).
\end{equation*}
This contradicts the inequality \eqref{eq:nu-ineq}.
Therefore, each cycle in the cycle decomposition of $g'$ contains
at most one element in $\{1,2,\dots,k\}$.
Namely, $g'$ is of the form
\begin{equation*}
g'=(\vj_k,k)\cdot\dots\cdot(\vj_{2},2)\cdot(\vj_1,1)\cdot h
\end{equation*}
for certain (possibly empty) disjoint strings
$\vj_1,\dots,\vj_k$ in $\{k+1,\dots,n\}$
and $h\in\sym n(\{k+1,\dots,n\})$.

For distinct elements $i_1,\dots,i_l\in\{1,2,\dots,k\}$,
we have
\begin{equation*}
(\vj_{i_l},i_l)\cdot\dotsb\cdot(\vj_{i_2},i_2)\cdot(\vj_{i_1},i_1)
\cdot(i_l,\dots,i_2,i_1)
=(\vj_{i_l},i_l,\dots,\vj_{i_2},i_2,\vj_{i_1},i_1).
\end{equation*}
This implies that $l$ distinct cycles in $g'$ turn into one cycle
in $g'\cdot(i_l,\dots,i_2,i_1)$, that is,
\begin{equation*}
\nu_n(g')-\nu_n(g'\cdot(i_l,\dots,i_2,i_1))=l-1.
\end{equation*}
Hence,
if $w\in\sym k$ is of the type $1^{r_1}2^{r_2}\dotsb k^{r_k}$,
then we have
\begin{equation*}
\nu_n(g')-\nu_n(g'w)=\sum_{l=1}^k r_l(l-1)=k-\nu_k(w).
\end{equation*}
Therefore it follows that
\begin{equation*}
\begin{split}
\sum_{w\in\sym k}\alpha^{n-\nu_n(g'w)}
&=\alpha^{n-\nu_n(g')}\sum_{w\in\sym k}\alpha^{k-\nu_k(w)}
=\alpha^{n-\nu_n(gw_0)}\prod_{1\le i<k}(1+i\alpha).
\end{split}
\end{equation*}
This completes the proof.
\end{proof}

Let us define the left action of $\sym m$
(resp. the right action of $\sym n$)
on the set $\Mat_{m,n}$
as permutations of row (resp. column) vectors:
\begin{equation*}
\begin{split}
\sigma\cdot(a_{ij})_{\substack{1\le i\le m\\ 1\le j\le n}}
&\deq(a_{\sigma^{-1}(i)j})_{\substack{1\le i\le m\\ 1\le j\le n}}
\qquad(\sigma\in\sym m),\\
(a_{ij})_{\substack{1\le i\le m\\ 1\le j\le n}}\cdot \tau
&\deq(a_{i \tau(j)})_{\substack{1\le i\le m\\ 1\le j\le n}}
\qquad(\tau\in\sym n).
\end{split}
\end{equation*}
Notice that $\sigma\cdot A=\pmat\sigma A$ and $A\cdot\tau=A\pmat\tau$
for $\sigma\in\sym m$, $\tau\in\sym n$ and $A\in\Mat_{m,n}$.
If $m=n$,
then we have
\begin{equation*}
\begin{split}
\adet(w\cdot A)&=\adet(a_{w^{-1}(i)j})
=\sum_{g\in\sym n}\alpha^{n-\nu_n(g)}\prod_{i=1}^n a_{w^{-1}g(i)i}\\
&=\sum_{g\in\sym n}\alpha^{n-\nu_n(wgw^{-1})}
\prod_{i=1}^n a_{g(i)w(i)}
=\adet(a_{i w(j)})=\adet(A\cdot w)
\end{split}
\end{equation*}
for any $w\in\sym n$ and any $A=(a_{ij})\in\Mat_n$.

\begin{lem}
The equality
\begin{equation*}
\sum_{w\in\sym n(I)}\adet(A\cdot w)=
\prod_{1\le i<k}(1+i\alpha)
\sum_{g\in\sym n}
\alpha^{m(g,I)}\prod_{i=1}^n a_{g(i)i}
\end{equation*}
holds for
$A=(a_{ij})_{1\le i,j\le n}\in\Mat_n$
and
$I\subset[n]$ such that $\card{I}=k$.
\end{lem}

\begin{proof}
Using Lemma \ref{lem:shifted_cycle_sum}, we have
\begin{equation*}
\begin{split}
\sum_{w\in\sym n(I)}\adet(A\cdot w)
&=\sum_{w\in\sym n(I)}\sum_{g\in\sym n}\alpha^{n-\nu_n(g)}
\prod_{i=1}^n a_{g(i)w(i)}
=\sum_{g\in\sym n}\sum_{w\in\sym n(I)}\alpha^{n-\nu_n(g)}
\prod_{i=1}^n a_{gw^{-1}(i)i}\\
&=\sum_{g\in\sym n}\ckakko{\sum_{w\in\sym n(I)}\alpha^{n-\nu_n(gw)}}
\prod_{i=1}^n a_{g(i)i}
=\prod_{1\le i<k}(1+i\alpha)\sum_{g\in\sym n}
\alpha^{m(g,I)}\prod_{i=1}^n a_{g(i)i}
\end{split}
\end{equation*}
as we desired.
\end{proof}

As a corollary, we have the following lemma.
\begin{lem}\label{lem:k-alternating}
For $I\subset[n]$ such that $\card{I}>k$
and $A\in\Mat_n$, the equalities
\begin{equation*}
\sum_{w\in\sym n(I)}\kdet(A\cdot w)=
\sum_{w\in\sym n(I)}\kdet(w\cdot A)=0
\end{equation*}
hold.
In particular, if $k+1$ column \pfkakko{row} vectors in $A$ are equal,
then $\kdet A=0$.
\qed
\end{lem}

Lemma \ref{lem:k-alternating} and the multilinearity of $\kdet$
yield immediately the
\begin{lem}\label{lem:key_lemma}
Let $A=(\va_1,\dots,\va_n)\in\Mat_n$.
If $\va_{i_1}=\dots=\va_{i_k}=\vb$ for some $1\le i_1<\dots<i_k\le n$,
then
\begin{equation*}
\kdet(\va_1,\dots,\va_j+\vb,\dots,\va_n)
=\kdet(\va_1,\dots,\va_j,\dots,\va_n)
\end{equation*}
for any $j\in[n]\setminus\{i_1,\dots,i_k\}$.
\qed
\end{lem}

When we regard $\sigma\in\sym n$ as an element in $\sym{n+m}$ ($m\in\N$) in natural way,
we notice that $\nu_{n+m}(\sigma)=\nu_n(\sigma)+m$.
Further, if we take a permutation $\tau\in\sym{m}$ and regard $\tau$
as an element in $\sym{n+m}$ which leave each letter in $[n]$ invariant,
then $\nu_{n+m}(\sigma\tau)=\nu_n(\sigma)+\nu_{m}(\tau)$.
This fact readily implies the following simple consequence
which will be used in the proof of Lemma \ref{lem:double_multiplex}
(see also Appendix).

\begin{lem}\label{lem:blockwise_multiplicativity}
The equality
\begin{equation*}
\adet\begin{pmatrix}
A_{11} & A_{12} \\ O & A_{22}
\end{pmatrix}=\adet(A_{11})\adet(A_{22})
\end{equation*}
holds.
In particular, $\adet(A_{11}\oplus A_{22})=\adet(A_{11})\adet(A_{22})$.
\end{lem}

\begin{proof}
Suppose that $A=(a_{ij})\in\Mat_{n+m}$ and
$A_{11}=(a_{ij})_{1\le i,j\le n}$,
$A_{22}=(a_{ij})_{n+1\le i,j\le n+m}$.
We also assume that $a_{ij}=0$ if $n+1\le i\le n+m$ and $1\le j\le n$.
Then it follows that
\begin{equation*}
\begin{split}
\adet A&=\sum_{\sigma\in\sym{n+m}}\alpha^{n+m-\nu_{n+m}(\sigma)}
\prod_{i=1}^{n+m}a_{i\sigma(i)}\\
&=\multsum{\sigma\in\sym{n+m}([n])\\ \tau\in\sym{n+m}(n+[m])}
\alpha^{n+m-\nu_{n+m}(\sigma\tau)}
\prod_{i=1}^{n}a_{i\sigma(i)}\prod_{i=1}^{m}a_{n+i,\tau(n+i)}\\
&=\multsum{\sigma\in\sym n\\ \tau\in\sym m}
\alpha^{n+m-\nu_{n}(\sigma)-\nu_m(\tau)}
\prod_{i=1}^{n}a_{i\sigma(i)}\prod_{i=1}^{m}a_{n+i,n+\tau(i)}
=\adet(A_{11})\adet(A_{22}).
\end{split}
\end{equation*}
This proves the claim.
\end{proof}

\section{Characterization of $-\frac1k$-determinants}
\label{sec:characterization}

In Lemma \ref{lem:k-alternating},
we prove that $\kdet$ has an alternating property
among $k+1$ column (and/or row) vectors.
In this section,
we show, conversely, this property
essentially characterizes $\kdet$.

We denote by $\A(\Mat_n(\C))$ the commutative $\C$-algebra
consisting of polynomial functions on $\Mat_n(\C)$.
The Lie algebra of $GL_n$ is denoted by $\liegl_n$,
and its universal enveloping algebra is denoted by $\U(\liegl_n)$.
The algebra $\A(\Mat_n(\C))$ has a $\U(\liegl_n)\times\sym n$-module structure by
defining
\begin{equation*}
(E_{ij}\cdot f)(X)=\sum_{k=1}^n x_{ik}\frac{\partial f}{\partial x_{jk}}(X)
\quad(1\le i,j\le n),\qquad
(\sigma\cdot f)(X)=f(X\cdot\sigma)\quad(\sigma\in\sym n)
\end{equation*}
for $f\in\A(\Mat_n(\C))$
where $E_{ij}$ are the standard basis of $\liegl_n$
and $x_{ij}$ are the standard coordinate functions on $\Mat_n(\C)$.
We note that this action of $\U(\liegl_n)$ is obtained
as the differential representation of $GL_n$ given by
$(g\cdot f)(X)=f(\tp{g}X)$ for $g\in GL_n$,
which is the contragradient representation
of the left regular representation on $\A(\Mat_n(\C))$.
Here $\tp{g}$ denotes the transposed matrix of $g$.

Let $\ML_n$ be a subspace of $\A(\Mat_n(\C))$
consisting of functions which are multilinear
with respect to \emph{column} vectors.
Clearly, we have
\begin{equation*}
\ML_n=\bigoplus_{1\le i_1,\dots,i_n\le n}\C\cdot x_{i_11}\dotsb x_{i_nn}.
\end{equation*}
The subspace $\ML_n$ is a $\U(\liegl_n)\times\sym n$-submodule of $\A(\Mat_n(\C))$.
For each $k\in\N$, we put
\begin{equation*}
\begin{split}
\AL_n^k\deq\set[4]{f\in\ML_n}
{I\subset[n],\,\card I>k
\then
\sum_{\tau\in\sym n(I)}f(X\cdot\tau)=0}
\end{split}
\end{equation*}
where $X=(x_{ij})_{1\le i,j\le n}$.
This subspace $\AL_n^k$ is also $\U(\liegl_n)$-invariant
because the actions of $\U(\liegl_n)$ and $\sym n$ on $\A(\Mat_n(\C))$
commutes each other.
We also see that $\AL_n^k$ is $\sym n$-invariant since
\begin{equation*}
\sum_{\tau\in\sym n(I)}(\sigma\cdot f)(X\cdot\tau)
=\sum_{\tau\in\sym n(I)}f(X\cdot\tau\sigma)
=\sum_{\tau\in\sym n(\sigma^{-1}I)}f(Y\cdot\tau)\Big|_{Y=X\cdot\sigma}=0
\end{equation*}
for any $I\subset[n]$, $\card I>k$
if $f\in\AL_n^k$ and $\sigma\in\sym n$.
Since $\kdet\in\AL_n^k$ by Lemma \ref{lem:k-alternating},
it follows that $\AL_n^k\supset\U(\liegl_n)\cdot\kdet(X)$.

\begin{thm}
The equality
$\AL_n^k=\U(\liegl_n)\cdot\kdet(X)$
holds for $k=1,2,\dots,n-1$.
\end{thm}

\begin{proof}
In \cite{MW2005}, it is shown that
\begin{equation}\label{eq:degeneration_of_cyclic_module}
\U(\liegl_n)\cdot\kdet(X)\cong
\bigoplus_{\substack{\lambda\vdash n\\ \lambda_1\le k}}
(\GLmod n\lambda)^{\oplus f^\lambda},
\end{equation}
where $\GLmod n\lambda$ denotes the highest weight $\U(\liegl_n)$-module of
highest weight $\lambda$,
which is the differential representation of $\GLmod n\lambda$
and we use the same symbol to indicate it.
The irreducible module $\GLmod n\lambda$ is realized
in $\U(\liegl_n)\cdot\kdet(X)$ as an image of
the Young symmetrizer
\begin{equation*}
c_T=\multsum{q\in C(T)\\ p\in R(T)}
\sgn(q)qp\in\C[\sym n]
\qquad(T\in\STab(\lambda)).
\end{equation*}
Here $C(T)$ and $R(T)$ are the column group and row group of $T$ respectively
(see, e.g. \cite{Weyl}).
Hence,
to prove the opposite inclusion $\AL_n^k\subset\U(\liegl_n)\cdot\kdet(X)$,
it is enough to show that each element $f$ in $\AL_n^k$
is killed by the Young symmetrizer $c_T$
when $T\in\STab(\lambda)$ and $\lambda_1>k$.
We now prove this.
The image $c_T\cdot f$ of $f\in\AL_n^k$ by $c_T$ is calculated as
\begin{equation*}
\begin{split}
(c_T\cdot f)(X)
&=\sum_{q\in C(T)}\sgn(q)\sum_{p\in R(T)}
f(X\cdot{qpq^{-1}q})
=\sum_{q\in C(T)}\sgn(q)\sum_{p\in R(qT)}
f(X\cdot{pq}).
\end{split}
\end{equation*}
For each $q\in C(T)$, we see that
\begin{equation*}
\sum_{p\in R(qT)}
f(X\cdot{pq})=
\sum_{p'\in R_1'(qT)}\ckakko{\sum_{p\in R_1(qT)}
(p'q\cdot f)(X\cdot{p})}=0
\end{equation*}
since $p'q\cdot f\in\AL_n^k$ by
$\sym n$-invariance of $\AL_n^k$.
Here $R_1(qT)$ is the subgroup of $R(qT)$ consisting of
permutations which moves only the entries in the first row of $qT$,
and $R_1'(qT)$ is the subgroup of $R(qT)$
which leave the first row of $qT$ invariant
so that $R(qT)=R_1(qT)\times R_1'(qT)$.
This completes the proof.
\end{proof}

\section{Determinants from a variation on wreath product groups}
\label{sec:multiple}

Let $m, n, k\in\N$.
For a matrix $A=(\va_1,\dots,\va_n)\in\Mat_{m,n}$,
we define the column $k$-plexing $A^{[k]}\in\Mat_{m,kn}$ of $A$ by
\begin{equation*}
\kple{A}\deq(\overbrace{\va_1,\dots,\va_1}^k,\dots,
\overbrace{\va_n,\dots,\va_n}^k).
\end{equation*}
This is nothing but the Kronecker product matrix
$A\otimes(1,\dots,1)$
of $A$ and $(1,\dots,1)\in\Mat_{1,k}$.
The row $k$-plexing $\kplerow{A}\in\Mat_{km,n}$ of $A$ is also defined
in a similar way.

\begin{ex}
If $A=
\begin{pmatrix}
a_1 & b_1 \\ a_2 & b_2 \\ a_3 & b_3
\end{pmatrix}\in\Mat_{3,2}$,
then
\begin{equation*}
\kple[2]A=\begin{pmatrix}
a_1 & a_1 & b_1 & b_1 \\
a_2 & a_2 & b_2 & b_2 \\
a_3 & a_3 & b_3 & b_3
\end{pmatrix}\in\Mat_{3,4},
\quad
\kple[3]A=\begin{pmatrix}
a_1 & a_1 & a_1 & b_1 & b_1 & b_1 \\
a_2 & a_2 & a_2 & b_2 & b_2 & b_2 \\
a_3 & a_3 & a_3 & b_3 & b_3 & b_3
\end{pmatrix}\in\Mat_{3,6}.
\end{equation*}
\end{ex}

We notice that
\begin{equation*}
\kple A=A\cdot\kple{(I_n)},\quad
\kplerow A=\kplerow{(I_m)}\cdot A
\end{equation*}
for $A\in\Mat_{m,n}$.
Hence one has the
\begin{lem}\label{lem:multiplex_relation}
Let $A\in\Mat_{m,n}$.
Then the equalities
\begin{equation*}
\kple{(PA)}=P\cdot\kple{A},\quad
\kplerow{(AQ)}=\kplerow{A}\cdot Q
\end{equation*}
hold for $P\in\Mat_m$, $Q\in\Mat_n$.
In particular, we have
\begin{equation*}
\sigma\cdot\kple A=\kple{\kakko{\sigma\cdot A}},\quad
\kplerow A\cdot \tau=\kplerow{\kakko{A\cdot \tau}}
\end{equation*}
for $\sigma\in\sym m$, $\tau\in\sym n$.
\qed
\end{lem}

\begin{dfn}
For a rectangular matrix
$A=(a_{ij})_{\substack{1\le i\le kn\\ 1\le j\le n}}\in\Mat_{kn,n}$,
we define the $k$-th \emph{wreath determinant} of $A$ by
\begin{equation*}
\wrdet A\deq\kdet(\kple A)
=\sum_{\sigma\in\sym{kn}}\kakko{-\frac1k}^{kn-\nu_{kn}(\sigma)}
\prod_{p=1}^n\prod_{l=1}^k a_{\sigma((p-1)k+l),p}.
\end{equation*}
\end{dfn}

By Lemma \ref{lem:key_lemma},
it is immediate to see that the equalities
\begin{equation*}
\begin{split}
\wrdet(\va_1,\dots,\va_{i-1},\va_i+c\va_j,\va_{i+1},\dots,\va_n)
&=\wrdet(\va_1,\dots,\va_{i-1},\va_i,\va_{i+1},\dots,\va_n)\quad(i\ne j),\\
\wrdet(\va_1,\dots,c\va_i,\dots,\va_n)
&=c^k\wrdet(\va_1,\dots,\va_i,\dots,\va_n)
\end{split}
\end{equation*}
hold for $A=(\va_1,\dots,\va_n)\in\Mat_{kn,n}$ and $c\in\C$.
Then it also follows that
\begin{equation}\label{eq:blockwise_alt}
\wrdet(A\cdot\sigma)
=(\sgn\sigma)^k\wrdet A
\qquad\kakko{\sigma\in\sym n}.
\end{equation}
In general, we have the
\begin{lem}\label{lem:key_lemma_2}
If $A\in\Mat_{kn,n}$ and $P\in \Mat_n$, then
\begin{equation*}
\wrdet(AP)=(\det P)^k \wrdet(A).
\end{equation*}
Namely, $\wrdet$ is a relative invariant of $GL_n$ in $\A(\Mat_{kn,n}(\C))$
with respect to the \pfkakko{right} regular representation
\pfkakko{See also Section \pfref{sec:detexp}}.
\qed
\end{lem}

\begin{ex}\label{ex:k-ply_multiplicativity}
Lemma \ref{lem:key_lemma_2} says that
the equality
\begin{align}\label{eq:Lemma_4.4_in_general}
\adet(\kple{(AP)})=(\det P)^k \adet(\kple A)
\end{align}
holds when $\alpha=-1/k$.
When $k=1$ and $\alpha=-1$,
this is nothing but the multiplicativity of the ordinary determinant.
We also notice that
\eqref{eq:Lemma_4.4_in_general} becomes trivial
when $\alpha=-1,-1/2,\dots,-1/(k-1)$.
Actually, because of Lemma \ref{lem:k-alternating},
each side of \eqref{eq:Lemma_4.4_in_general} vanishes
for such values.
Further, we notice that
\eqref{eq:Lemma_4.4_in_general} holds
only if $\alpha=-1,-1/2,\dots,-1/k$.
Actually, if $\adet(\kple X)$ satisfies
\eqref{eq:Lemma_4.4_in_general},
then the ratio $\adet(\kple X)/\wrdet(X)$ gives
an absolute invariant of $GL_n$,
which must be a constant.
If the constant is $0$,
then it follows from \eqref{eq:adet_of_1} that
$\alpha=-1,-1/2,\dots,-1/(k-1)$.
If the constant is not $0$,
then we immediately have $\alpha=-1/k$.
Here we give a simple and direct example.
When $n=k=2$ and $P=\left(\begin{smallmatrix}
1 & 1\\ 0 & 1\end{smallmatrix}\right)$, we have
\begin{equation*}
\begin{split}
&\adet(\kple[2]{(AP)})-(\det P)^2 \adet(\kple[2] A)\\
=&(1+\alpha)(1+2\alpha)
\kakko{(1+3\alpha)a_{11}a_{21}a_{31}a_{41}%
+2\alpha(a_{12}a_{21}+a_{11}a_{22})a_{31}a_{41}%
+(1+\alpha)a_{11}a_{21}(a_{32}a_{41}+a_{31}a_{42})}
\end{split}
\end{equation*}
which is identically zero only if $\alpha=-1,-\frac12$.
See also Corollary \ref{cor:char_of_wrdet}.
\end{ex}

\begin{lem}\label{lem:double_multiplex}
If $A\in\Mat_n$, then the equality
\begin{equation*}
\kdet\kakko{\kplerc{A}}
=\wrdet\kakko{\kplerow{A}}
=\kakko{\frac{k!}{k^k}}^n(\det A)^k
\end{equation*}
holds for any $k\in\N$.
\end{lem}

\begin{proof}
By Lemmas \ref{lem:multiplex_relation} and \ref{lem:key_lemma_2},
we have
\begin{equation*}
\begin{split}
\kdet\kakko{\kplerc{A}}
&=\wrdet\kakko{\kplerow{A}}
=\wrdet\kakko{\kplerow{(I_n)}\cdot A}\\
&=\wrdet\kakko{\kplerow{(I_n)}}\cdot(\det A)^k
=\kdet\kakko{\kplerc{(I_n)}}\cdot(\det A)^k.
\end{split}
\end{equation*}
Since $\kplerc{(I_n)}=\overbrace{\1_k\oplus\dots\oplus\1_k}^n$
and $\kdet(\1_k)=\prod_{1\le i<k}(1-\frac ik)=\frac{k!}{k^k}$,
we have $\kdet\kakko{\kplerc{(I_n)}}=\kakko{\frac{k!}{k^k}}^n$
by Lemma \ref{lem:blockwise_multiplicativity}.
This completes the proof.
\end{proof}
This Lemma will be used in \S7.

We consider the two injective homomorphisms
$\phi:\sym k^n\to\sym{kn}$ and $\psi:\sym n\to\sym{kn}$ defined as 
\begin{equation*}
\begin{split}
\phi(\sigma_1,\dots,\sigma_n)&:
[kn]\ni(i-1)k+j\longmapsto(i-1)k+\sigma_i(j)\in[kn]
\quad(1\le i\le n,\,1\le j\le k),\\
\psi(\tau)&:
[kn]\ni(i-1)k+j\longmapsto(\tau(i)-1)k+j\in[kn]
\quad(1\le i\le n,\,1\le j\le k)
\end{split}
\end{equation*}
for $(\sigma_1,\dots,\sigma_n)\in\sym k^n$ and $\tau\in\sym n$.
To avoid the confusion, we put
$S_k^n\deq\phi(\sym k^n)$ and $S_n\deq\psi(\sym n)$.
We note that $S_k^n$ is the Young subgroup $\sym{(k^n)}$
of $\sym{kn}$ corresponding to the partition $(k^n)\vdash kn$.

By the definition of $k$-plexing, one finds that 
$\kple A\cdot\sigma=\kple A$ for $A\in\Mat_{kn,n}$
and $\sigma\in S_k^n$,
whence it follows that
\begin{equation*}
\wrdet\kakko{\sigma\cdot A}
=\kdet\kakko{\sigma\cdot\kple A}
=\kdet\kakko{\kple A\cdot\sigma}
=\kdet\kakko{\kple A}
=\wrdet A
\qquad\kakko{\sigma\in S_k^n}.
\end{equation*}
We also see that
$\kple A\cdot\psi(\tau)=\kple{(A\cdot\tau)}$
for any $\tau\in\sym n$.
Hence we have 
\begin{equation*}
\begin{split}
\wrdet\kakko{\psi(\tau)\cdot A}
&=\kdet\kakko{\psi(\tau)\cdot\kple A}
=\kdet\kakko{\kple A\cdot\psi(\tau)}\\
&=\wrdet\kakko{A\cdot \tau}
=(\sgn \tau)^k\wrdet A
\qquad\kakko{\tau\in\sym n}
\end{split}
\end{equation*}
by \eqref{eq:blockwise_alt}.
Consequently, we obtain the
\begin{lem}\label{lem:relative_invariance_of_wrdet}
If $A\in\Mat_{kn,n}$,
then
\begin{equation*}
\wrdet\kakko{g\cdot A}
=\chi_{n,k}(g)^k\wrdet A
\end{equation*}
for any $g\in\wsym kn$.
In other words,
$\C\cdot\wrdet\subset\A(\Mat_{kn,n})$
defines a one-dimensional representation
of $\wsym kn$.
Here $\wsym kn\deq S_k^n\rtimes S_n$ is the wreath product group
\pfkakko{see \pfcite{Mac}}.
The character $\chi_{n,k}$ of $\wsym kn$ is defined by
\begin{equation*}
\chi_{n,k}(g)=\sgn \tau
\end{equation*}
for $g=(\phi(\sigma_1,\dots,\sigma_n);\psi(\tau))$
\pfkakko{$\sigma_i\in\sym k,\,\tau\in\sym n$}.
\qed
\end{lem}

\section{Expressions of wreath determinants and
$(GL_{kn},\, GL_n)$-duality}
\label{sec:detexp}

For given two linear spaces $V$ and $W$,
as a $GL(V)\times GL(W)$-module,
the multiplicity-free decomposition
\begin{equation}\label{eq:(GL,GL)-duality}
\symp(V\otimes W)\cong\bigoplus_{\lambda}\GLmod V\lambda\boxtimes\GLmod W\lambda
\end{equation}
of the symmetric algebra $\symp(V\otimes W)$ holds. 
Here $\lambda$ runs over the partitions
such that $\len\lambda\le\min\{\dim V,\dim W\}$.
This fact is referred as \emph{$(GL(V),GL(W))$-duality}
(see \cite{Howe1} and \cite{Weyl}).

The algebra $\A(\Mat_{kn,n})$ has a
$GL_{kn}\times GL_n$-module structure given by
\begin{equation*}
((g,h).f)(A)\deq f(\tp{g}Ah)
\qquad(g\in GL_{kn},\, h\in GL_n,\,A\in\Mat_{kn,n}),
\end{equation*}
where $\tp{g}$ denotes the transposition of $g$
with respect to the standard coordinate.
We see that
\begin{equation*}
\A(\Mat_{kn,n})\cong\A\kakko{(\C^{kn})^*\otimes(\C^{n})^*}
\cong\symp\kakko{\C^{kn}\otimes\C^{n}}
\end{equation*}
as $GL_{kn}\times GL_n$-module.
Here $V^*$ indicates the
contragradient representation of $V$.
We notice that if $(\rho, V)$ is a representation of $GL_m$,
then $\tilde\rho(g)=\rho(\tp{g}^{-1})$ ($g\in GL_m$)
defines a representation on $V$
which is equivalent to $V^*$.

\begin{rem}
It is standard to define a representation of $GL_{kn}\times GL_n$
on the algebra $\A(\Mat_{kn,n}(\C))$ by
\begin{equation*}
((g,h).f)(A)\deq f(g^{-1}Ah)
\qquad(g\in GL_{kn},\, h\in GL_n,\,A\in\Mat_{kn,n}),
\end{equation*}
which is a combination of the left regular action of $GL_{kn}$
and the right regular action of $GL_n$.
If we adopt this one, however,
then it is no longer a \emph{polynomial} representation.
Instead, in our argument, we adopt the contragradient of
the left regular action of $GL_{kn}$ 
so that each (irreducible) factor of
the $GL_{kn}\times GL_n$-module $\A(\Mat_{kn,n}(\C))$
is polynomial.
\end{rem}

By $(GL_{kn},GL_n)$-duality,
one has the multiplicity-free decomposition of $\A(\Mat_{kn,n})$: 
\begin{equation*}
\A(\Mat_{kn,n})\cong\bigoplus_{\len\lambda\le n}
\GLmod{kn}\lambda\boxtimes \GLmod n\lambda.
\end{equation*}
If we look at the $\det$-eigenspace with respect to
the left action of the diagonal torus $T_{kn}\cong(\C^{\times})^{kn}$ of $GL_{kn}$,
then we have
\begin{equation*}
\A(\Mat_{kn,n})^{T_{kn},\det}\cong\bigoplus_{\len\lambda\le n}
\kakko{\GLmod{kn}\lambda}^{T_{kn},\det}\boxtimes \GLmod n\lambda.
\end{equation*}
Here, for a $GL_{kn}$-module $V$,
we denote by $V^{T_{kn},\det}$
the $\det$-eigenspace
\begin{equation*}
V^{T_{kn},\det}=\set[2]{v\in V}{t.v=\det(t)v\quad(t\in T_{kn})}
\end{equation*}
with respect to $T_{kn}$.
Since the symmetric group $\sym{kn}$ is the normalizer of $T_{kn}$ in $GL_{kn}$,
each $\det$-eigenspace $\kakko{\GLmod{kn}\lambda}^{T_{kn},\det}$
becomes a $\sym{kn}$-module.
It is known that the equivalence
$\kakko{\GLmod{kn}\lambda}^{T_{kn},\det}\cong\Smod{kn}\lambda$
holds as $\sym{kn}$-modules
if $\lambda$ is a partition of $kn$
(see, e.g. \cite{Howe1}).

Let us denote by $M_{n,k}$ the irreducible $GL_{kn}\times GL_n$-submodule
of $\A(\Mat_{kn,n})$
corresponding to the partition $(k^n)$, that is,
$M_{n,k}\cong \GLmod{kn}{(k^n)}\boxtimes \GLmod n{(k^n)}$.
As $\sym{kn}$-modules,
we have the equivalence
\begin{equation*}
M_{n,k}^{T_{kn},\det}
\cong \kakko{\GLmod{kn}{(k^n)}}^{T_{kn},\det}
\boxtimes \GLmod n{(k^n)}
\cong \kakko{\GLmod{kn}{(k^n)}}^{T_{kn},\det}
\cong \Smod{kn}{(k^n)}
\end{equation*}
since the multiplicity space $\GLmod n{(k^n)}$ is
of dimension one.
In particular, we have
$\dim M_{n,k}^{T_{kn},\det}=f^{(k^n)}$.

By Lemma \ref{lem:key_lemma_2} and
$(GL_{kn},GL_n)$-duality,
it follows that $\wrdet\in M_{n,k}$.
Moreover, since
\begin{equation*}
(\diag(c_1,\dots,c_{kn}).\wrdet)(A)=
\wrdet(\tp\diag(c_1,\dots,c_{kn})A)=
\kakko{\prod_{i=1}^{kn}c_i} \wrdet A,
\end{equation*}
it follows that $\wrdet$ belongs to
$M_{n,k}^{T_{kn},\det}$.

For each standard tableau
$T=(t_{ij})_{\substack{1\le i\le n\\ 1\le j\le k}}\in\STab((k^n))$,
we define the function $\tdet_T$ on $\Mat_{kn,n}$ by
\begin{equation*}
\tdet_T(A)\deq\prod_{l=1}^k \det(a_{t_{il},j})_{1\le i,j\le n}
\qquad(A=(a_{ij})_{\substack{1\le i\le kn\\ 1\le j\le n}}\in\Mat_{kn,n}).
\end{equation*}
We also define the matrix $I(T)\in\Mat_{kn,n}$
so that $t_{ij}$-th row vector of $I(T)$ is equal to the $i$-th
fundamental row vector $\ve_i=(0,\dots,0,\overset{\text{$i$-th}}1,0,\dots,0)$
for each $i=1,\dots,n$ and $j=1,\dots,k$.
In other words,
if we define $g(T)\in\sym{kn}$ for $T\in\STab((k^n))$ by
\begin{equation}\label{eq:def_of_g(T)}
g(T)\kakko[1]{(i-1)k+j}=t_{ij}\quad(1\le i\le n,\,1\le j\le k),
\end{equation}
then $I(T)=g(T)\cdot\kplerow{(I_n)}$.
Denote by $T_0$ the standard tableau with shape $(k^n)$
whose $(i,j)$-entry is $(i-1)k+j$.
We note that
$g(T)\in\sym{kn}$ is the permutation determined by $g(T)\cdot T_0=T$
for each $T\in\STab((k^n))$.

\begin{lem}\label{lem:tdet_formula}
For $T,U\in\STab((k^n))$, the equality
\begin{equation*}
\tdet_T(I(U))=\begin{cases}
1 & T=U\\
0 & T\ne U
\end{cases}
\end{equation*}
holds.
\end{lem}

\begin{proof}
When $T=U$,
the $t_{il}$-th row vector
$I(T)_{t_{il}}$ of $I(T)$ is equal to $\ve_i$
if $i\in[n]$ and $l\in[k]$,
and hence $\tdet_T(I(T))=1$.
When $T=(t_{ij})$ and $U=(u_{ij})$ are
distinct standard tableaux of shape $(k^n)$,
there exists a pair $(s_1,s_2)$ of distinct elements in $[kn]$ such that
$s_1$ and $s_2$ are in the same column of $T$
and in the same row of $U$, say
$s_1=t_{i_1c}=u_{rj_1}$ and $s_2=t_{i_2c}=u_{rj_2}$
($i_1\ne i_2$, $j_1\ne j_2$).
Then we have
\begin{equation*}
I(U)_{t_{i_1}}=I(U)_{t_{i_2}}=\ve_r,
\end{equation*}
which implies that $\det(I(U)_{t_{ic},j})_{1\le i,j\le n}=0$,
and hence $\tdet_T(I(U))=0$.
\end{proof}

\begin{thm}\label{thm:detexp_of_wrdet}
The wreath determinant $\wrdet A$ of a matrix $A\in\Mat_{kn,n}$
is expressed as a linear combination
\begin{equation*}
\wrdet A=\sum_{T\in\STab((k^n))} \wrdet I(T)\cdot \tdet_T(A)
\end{equation*}
of $\tdet_T(A)$ for $T\in \STab((k^n))$.
The coefficient $\wrdet I(T)$ is given by the sum
\begin{equation*}
\wrdet I(T)=\sum_{\sigma\in S_k^n}\kakko{-\frac1k}^{kn-\nu_{kn}(g(T)\sigma)},
\end{equation*}
where $g(T)\in\sym{kn}$ is a permutation
defined by \eqref{eq:def_of_g(T)}.
\end{thm}

\begin{proof}
We observe that
$\tdet_T(A)$ is a homogeneous polynomial in $a_{ij}$ of degree $kn$
satisfying the condition that $\tdet_T(AP)=(\det P)^k\tdet_T(A)$
for any $P\in\Mat_n$.
We also see that
\begin{equation*}
(\diag(c_1,\dots,c_{kn}).\tdet_T)(A)=
\tdet_T(\tp\diag(c_1,\dots,c_{kn})A)=
\kakko{\prod_{i=1}^{kn}c_i} \tdet_T A.
\end{equation*}
Thus, it follows that
every $\tdet_T$ belongs to $M_{n,k}^{T_{kn},\det}$
by $(GL_{kn},GL_n)$-duality.

We show that
$\{\tdet_T\}_{T\in\STab((k^n))}$ are linearly independent.
Suppose that
\begin{equation*}
\sum_{T\in\STab((k^n))}C_T\tdet_T(A)=0
\end{equation*}
for any $A\in\Mat_{kn,n}$.
Then, by Lemma \ref{lem:tdet_formula}, we have
\begin{equation*}
0=\sum_{T\in\STab((k^n))}C_T\tdet_T(I(U))=C_U
\end{equation*}
for each $U\in\STab((k^n))$,
which assures the linear independence of
$\{\tdet_T\}_{T\in\STab((k^n))}$.
Since $\dim M_{n,k}^{T_{kn},\det}=f^{(k^n)}$,
it follows that $\{\tdet_T\}_{T\in\STab((k^n))}$ is a
basis of $M_{n,k}^{T_{kn},\det}$.
Hence $\wrdet$ is written as
\begin{equation*}
\wrdet A=\sum_{T\in\STab((k^n))} C'_T \tdet_T(A)
\qquad(A\in\Mat_{kn,n}).
\end{equation*}
By Lemma \ref{lem:tdet_formula} again,
the coefficient $C'_U$ for $U\in\STab((k^n))$ is
calculated as
\begin{equation*}
\wrdet I(U)=\sum_T C'_T\tdet_T(I(U))=C'_U \tdet_U(I(U))=C'_U.
\end{equation*}
This completes the proof of the theorem.
(The coefficient $\wrdet I(U)$ is calculated later in
Section \ref{sec:sgn}.)
\end{proof}

\begin{ex}\label{ex:detexp_of_wrdet_for_(n,k)=(3,2)}
When $n=3$ and $k=2$, there are five standard tableaux with shape $(2^3)$:
\begin{equation*}
U_1=\young(12,34,56)\,,\quad
U_2=\young(12,35,46)\,,\quad
U_3=\young(13,24,56)\,,\quad
U_4=\young(13,25,46)\,,\quad
U_5=\young(14,25,36)\,.
\end{equation*}
(We remark that $T_0=U_1$ in this case.)
The corresponding matrices $I(U_p)$ are given by
\begin{equation*}
\begin{split}
&I(U_1)=\begin{pmatrix}
1 & 0 & 0 \\
1 & 0 & 0 \\
0 & 1 & 0 \\
0 & 1 & 0 \\
0 & 0 & 1 \\
0 & 0 & 1
\end{pmatrix},\quad
I(U_2)=\begin{pmatrix}
1 & 0 & 0 \\
1 & 0 & 0 \\
0 & 1 & 0 \\
0 & 0 & 1 \\
0 & 1 & 0 \\
0 & 0 & 1
\end{pmatrix},\quad
I(U_3)=\begin{pmatrix}
1 & 0 & 0 \\
0 & 1 & 0 \\
1 & 0 & 0 \\
0 & 1 & 0 \\
0 & 0 & 1 \\
0 & 0 & 1
\end{pmatrix},\\
&I(U_4)=\begin{pmatrix}
1 & 0 & 0 \\
0 & 1 & 0 \\
1 & 0 & 0 \\
0 & 0 & 1 \\
0 & 1 & 0 \\
0 & 0 & 1
\end{pmatrix},\quad
I(U_5)=\begin{pmatrix}
1 & 0 & 0 \\
0 & 1 & 0 \\
0 & 0 & 1 \\
1 & 0 & 0 \\
0 & 1 & 0 \\
0 & 0 & 1
\end{pmatrix},
\end{split}
\end{equation*}
and their $2$-wreath determinants are calculated as
\begin{equation*}
\wrdet[2]I(U_1)=\frac18,\quad
\wrdet[2]I(U_2)=\wrdet[2]I(U_3)=-\frac1{16},\quad
\wrdet[2]I(U_4)=\wrdet[2]I(U_5)=\frac1{32}.
\end{equation*}
Thus we have
\begin{equation*}
\wrdet[2]A=\frac18\tdet_{U_1}(A)-\frac1{16}\tdet_{U_2}(A)
-\frac1{16}\tdet_{U_3}(A)
+\frac1{32}\tdet_{U_4}(A)+\frac1{32}\tdet_{U_5}(A)
\end{equation*}
for $A\in\Mat_{6,3}$.
\end{ex}

As a corollary of the theorem, we obviously have the
\begin{cor}
For $A\in\Mat_{p,n}$ and $B\in\Mat_{q,n}$,
we denote by $A\boxplus B\in\Mat_{p+q,n}$ the matrix
obtained by piling $A$ on $B$.
If $A_1,\dots,A_k\in\Mat_{n,n}$, then the equality
\begin{equation*}
\wrdet(A_1\boxplus\dotsb\boxplus A_k)
=\sum_{T\in\STab((k^n))}\wrdet I(T)\prod_{i=1}^k \det B_i(T)
\end{equation*}
holds, where $B_j(T)$ is a matrix whose $i$-th row is
equal to the $t_{ij}$-th row of $A_1\boxplus\dotsb\boxplus A_k$.
\end{cor}

\begin{ex}
If $A=\begin{pmatrix}
a_{11} & a_{12} \\ a_{21} & a_{22}
\end{pmatrix}$
and
$B=\begin{pmatrix}
b_{11} & b_{12} \\ b_{21} & b_{22}
\end{pmatrix}$,
then we have
\begin{equation*}
\begin{split}
\wrdet[2]\kakko{A\boxplus B}=
\wrdet[2]\begin{pmatrix}
a_{11} & a_{12} \\
a_{21} & a_{22} \\
b_{11} & b_{12} \\
b_{21} & b_{22}
\end{pmatrix}
&=\frac14
\begin{vmatrix}
a_{11} & a_{12} \\
a_{21} & a_{22}
\end{vmatrix}
\begin{vmatrix}
b_{11} & b_{12} \\
b_{21} & b_{22}
\end{vmatrix}
-\frac18
\begin{vmatrix}
a_{11} & a_{12} \\
b_{11} & b_{12}
\end{vmatrix}
\begin{vmatrix}
a_{21} & a_{22} \\
b_{21} & b_{22}
\end{vmatrix}.
\end{split}
\end{equation*}
\end{ex}

Recall that the wreath determinant $\wrdet$ is $S_k^n$-invariant.
By the Frobenius reciprocity,
it follows that
\begin{equation*}
\dim\kakko{M_{n,k}^{T_{kn},\det}}^{S_k^n}
=\inprod{\res_{\sym{kn}}^{S_k^n}\kakko{M_{n,k}^{T_{kn},\det}}}
{\trivrpn{S_k^n}}_{S_k^n}
=\inprod{M_{n,k}^{T_{kn},\det}}{\ind_{S_k^n}^{\sym{kn}}\trivrpn{S_k^n}}_{\sym{kn}}
=K_{(k^n)(k^n)}=1,
\end{equation*}
where $\inprod{V}{W}_{G}$ denotes the intertwining number
of two $G$-modules $V$ and $W$,
and $\trivrpn{G}$ is the trivial representation of $G$.
Hence we have
\begin{equation}\label{eq:S_k^n-invariant_space_is_one_dimensional}
\kakko{M_{n,k}^{T_{kn},\det}}^{S_k^n}=\C\cdot\wrdet(X).
\end{equation}
This fact implies that
$\sum_{\sigma\in S_k^n}f(\sigma\cdot X)$
is proportional to $\wrdet(X)$ for any $f\in M_{n,k}^{T_{kn},\det}$.
Therefore, we have
\begin{equation*}
\sum_{\sigma\in S_k^n}\tdet_{T_0}(\sigma\cdot X)
=C\wrdet(X)
\end{equation*}
for a certain constant $C$.
If we set $X=\kplerow{(I_n)}$,
then we have
\begin{equation*}
C=\frac1{\wrdet\kplerow{(I_n)}}
\sum_{\sigma\in S_k^n}\tdet_{T_0}\kakko{\sigma\cdot\kplerow{(I_n)}}
=\kakko{\frac{k^k}{k!}}^{\!n}\sum_{\sigma\in S_k^n}1
=k^{kn}.
\end{equation*}
Consequently, we obtain another (symmetric) expression of $\wrdet(X)$ as follows.
\begin{cor}\label{cor:symmetric_exp_of_wrdet}
The equality
\begin{equation*}
\wrdet(A)=\frac1{k^{kn}}\sum_{\sigma\in S_k^n}\tdet_{T_0}(\sigma\cdot A)
\end{equation*}
holds for any $A\in\Mat_{kn,n}$.
\qed
\end{cor}

As a corollary of the discussion above, we obtain the
\begin{cor}[Characterization of the wreath determinant]
\label{cor:char_of_wrdet}
Put
\begin{equation*}
\A(\Mat_{kn,n})^{\chi_{n,k}^k,\det^k\!}=
\set[2]{f\in\A(\Mat_{kn,n})}{f(\sigma XP)=\chi_{n,k}(\sigma)^k(\det P)^k f(X),\
\sigma\in\wsym kn,\, P\in GL_n}.
\end{equation*}
Then $\A(\Mat_{kn,n})^{\chi_{n,k}^k,\det^k\!}$ is a
one dimensional subspace spanned by $\wrdet$.
Namely, the equality
\begin{equation*}
\A(\Mat_{kn,n})^{\chi_{n,k}^k,\det^k\!}=\C\cdot\wrdet(X)
\end{equation*}
holds.
\qed
\end{cor}

Corollary \ref{cor:char_of_wrdet} and
Example \ref{ex:k-ply_multiplicativity}
suggest the following problem:
Describe the irreducible decomposition and singular values
of the cyclic module
$\U(\liegl_{kn})\cdot\adet(\kple X)\subset\A(\Mat_{kn,n})$
($X=(x_{ij})_{1\le i\le kn,1\le j\le n}$).
This is solved in the following way.
If $\alpha=0$, then we see that
\begin{align*}
\U(\liegl_{kn})\cdot\adet[0](\kple X)
\cong\symp^k(\C^{kn})^{\otimes n}
\cong\bigoplus_{\lambda\vdash kn}
\kakko{\GLmod{kn}\lambda}^{\oplus K_{\lambda,(k^n)}}
\end{align*}
by a similar discussion in \cite{KMW}
(we also refer to \cite{MW2005} for the case where $k=1$).
By \cite{MW2005},
the $\lambda$-isotypic component of the module
$\U(\liegl_{kn})\cdot\adet(\widetilde X)\subset\A(\Mat_{kn})$
does have a positive multiplicity
if and only if $f_\lambda(\alpha)\ne0$
and is given by $\U(\liegl_{kn})\cdot\Imm_\lambda(\widetilde X)$
(We put $\widetilde X=(x_{ij})_{1\le i,j\le kn}$
in order to avoid confusion).
Here $\Imm_\lambda(\widetilde X)$ is
the \emph{immanant} of $\widetilde X$
for $\lambda$ and 
$f_\lambda(\alpha)\deq\prod_{(i,j)\in\lambda}(1+(j-i)\alpha)$
is the (modified) content polynomial for $\lambda$.
Since the map
$\A(\Mat_{kn})\ni f(\widetilde X)
\mapsto f(\kple X)\in\A(\Mat_{kn,n})$
defines a $GL_{kn}$-intertwiner,
we see that
\begin{align*}
\text{the $\lambda$-isotypic component of
$\U(\liegl_{kn})\cdot\adet(\kple X)$}\cong
\begin{cases}
\U(\liegl_{kn})\cdot\Imm_\lambda(\kple X) & f_\lambda(\alpha)\ne0\\
0 & \text{otherwise}
\end{cases}
\end{align*}
for $\lambda\vdash kn$.
Thus it follows that
$\U(\liegl_{kn})\cdot\Imm_\lambda(\kple X)\cong
\kakko{\GLmod{kn}\lambda}^{\oplus K_{\lambda,(k^n)}}$.
Hence we obtain the following theorem
which is regarded as a generalization of the result in \cite{MW2005}.
\begin{thm}
The irreducible decomposition of the cyclic module
generated by $\adet(\kple X)$ is given by
\begin{align*}
\U(\liegl_{kn})\cdot\adet(\kple X)
\cong\bigoplus_{\substack{\lambda\vdash kn\\ f_\lambda(\alpha)\ne0}}
\kakko{\GLmod{kn}\lambda}^{\oplus K_{\lambda,(k^n)}}.
\end{align*}
In particular,
the singular values are given as roots of
the content polynomials.
\qed
\end{thm}

\subsection*{Remarks on this section}

Let $\symp(\C^n)=\sum_{k\geq 0}\symp^k(\C^n)$
be the homogeneous decomposition of $\symp(\C^n)$.
Each symmetric power $\symp^k(\C^n)$, that is,
the space of $k$-th symmetric tensors
defines an irreducible $GL_n(\C)$-module \cite{F}.
We see that the eigenspace decomposition of
the $GL_m\times GL_n$-module $\symp(\C^m\otimes\C^n)$
with respect to the diagonal torus $T_m$ of $GL_m(\C)$ is
given by
\begin{equation*}
\symp(\C^m\otimes\C^n)\cong\bigoplus_{k_1,\dots,k_m\ge0}
\symp^{k_1}(\C^n)\otimes\dotsb\otimes\symp^{k_m}(\C^n).
\end{equation*}
Hence the $m$-th tensor product
$\symp^k(\C^n)^{\otimes m}$ can be identified to
the $\det^k$-eigenspace
\begin{equation*}
\symp(\C^m\otimes \C^n)^{T_m,\det^k}
=\set[2]{v\in\symp(\C^m\otimes\C^n)}{t.v=(\det t)^kv\ (t\in T_m)}
\end{equation*}
for $T_m$ \cite{Howe1}.
By $(GL_m,\, GL_n)$-duality \eqref{eq:(GL,GL)-duality},
we see that
\begin{equation}\label{Wreath-GL-duality}
\symp^k(\C^n)^{\otimes m}\cong
\symp(\C^m\otimes \C^n)^{T_m,\det^k}\cong
\sum_{\len\lambda\leq \min\{m,n\}}
(\GLmod m\lambda)^{T_m,\det^k}\boxtimes \GLmod n\lambda.
\end{equation}
We notice that
$(\GLmod m\lambda)^{T_m,\det^k}=\{0\}$ unless $\lambda\vdash km$,
and hence the last sum \eqref{Wreath-GL-duality} is effectively
over the partitions of $km$.
Note also that
$\dim (\GLmod m\lambda)^{T_m,\det^k}=K_{\lambda (k^m)}$ and
$(\GLmod m\lambda)^{T_m,\det^k}$ is stable under
the action of the Weyl group $\sym m$ of $GL_m(\C)$.
We note that
the decomposition \eqref{Wreath-GL-duality}
for $k=1$ gives $(\sym m,\, GL_n)$-duality (Schur duality)
\begin{equation}\label{S-duality}
(\C^n)^{\otimes m} \cong \sum_{\lambda\vdash m,\, \ell(\lambda)\leq n}
\Smod m\lambda\boxtimes \GLmod n\lambda.
\end{equation}

Suppose now $\lambda\vdash km$.
The group $\wsym km$
acts on the weight space $(\GLmod m\lambda)^{T_m,\det^k}$
because
the wreath product $\wsym km=S_k^m\rtimes S_m$
is obviously acting on the space
$\symp^k(\C^n)^{\otimes m}$.
Since $\sym k$ acts on
$\symp^k(\C^n)^{\otimes m}$ trivially,
its action on the weight space $(\GLmod m\lambda)^{T_m,\det^k}$
is also trivial.
Hence, \eqref{Wreath-GL-duality} does not provide the
irreducible decomposition as a bi-module of
$(\wsym km,\, GL_n(\C))$. Then, the question
how the space $(\GLmod m\lambda)^{T_m,\det^k}$ decomposes
as a $\sym m$-module comes into being. Now we establish this question
in a concrete way.
From Schur duality, as a $\sym m\times GL(\symp^k(\C^n))$-module,
we obtain
\begin{equation*}
\symp^k(\C^n)^{\otimes m}
\cong \sum_{\mu\vdash m, \, \ell(\mu)\leq N}
\Smod m\mu\boxtimes \GLmod N\mu,
\end{equation*}
where $N=\dim \symp^k(\C^n)=\binom{n+k-1}{k}\geq n$.
Decompose the module $\GLmod N\lambda$ of $GL(\symp^k(\C^n))$
into irreducible ones
as a representation of the subgroup $GL_n(\C)$
of $GL(\symp^k(\C^n))$:
\begin{equation*}
\GLmod N\mu\Big|_{GL_n(\C)}\cong \sum_{\lambda,\,\ell(\lambda)\leq n}
\kakko{\GLmod n\lambda}^{\oplus m_\lambda(\mu)},
\end{equation*}
$m_\lambda(\mu)$ being the multiplicity of $\GLmod n\lambda$ in the
irreducible summands of the restriction.
Then we have
\begin{equation*}
\symp^k(\C^n)^{\otimes m}\cong
\sum_{\lambda,\,\ell(\lambda)\leq n} \;\sum_{\mu\vdash m}
\kakko{\Smod m\mu \boxtimes \GLmod n\lambda}^{\oplus m_\lambda(\mu)}.
\end{equation*}
Therefore, it follows from \eqref{Wreath-GL-duality} that
\begin{equation}\label{eq:decomp_of_J}
\sum_{\mu\vdash m} 
\kakko{\Smod m\mu}^{\oplus m_\lambda(\mu)}
\cong (\GLmod m\lambda)^{T_m,\det^k}.
\end{equation}
The procedure explained above is a special case of
the problem for computing \emph{plethysm}
(or the functorial composition of operations
$\lambda\mapsto \GLmod{}\lambda$) (see \cite{Mac}, \cite{Howe2}).
Note also that the problem for describing the decomposition
\eqref{eq:decomp_of_J} for $\lambda\vdash km$ explicitly
comes up naturally when one wants to know the
structure of the cyclic $GL_n(\C)$-module generated by
$\adet(X)^k\;(k=(1,)2,3,\ldots)$
(see \cite{KMW}).

\section{Formulas for wreath determinants \`a la Cauchy et van der Monde}
\label{sec:cauchy}

We give an analogue of the Cauchy determinant formula
in the context of wreath determinants developed in the previous sections.
\begin{prop}
Let $n, k\in\N$ and
$x_1,\dots,x_{kn},y_1,\dots,y_n$ be commutative variables.
Put
\begin{equation*}
\begin{split}
&C_{n,k}(x,y)=
\kakko{\frac1{x_i+y_j}}_{\substack{1\le i\le kn\\ 1\le j\le n}},\quad
V_{n,k}(x)=(x_i^{n-j})_{\substack{1\le i\le kn\\ 1\le j\le n}}.
\end{split}
\end{equation*}
Then we have
\begin{equation}\label{eq:Cauchy-Vandermonde}
\wrdet C_{n,k}(x,y)
=\frac{\Delta_n(y)^k}
{\displaystyle\multprod{1\le i\le kn\\1\le j\le n}(x_i+y_j)}
\wrdet V_{n,k}(x).
\end{equation}
Here $\Delta_n(y)$ denotes the difference product
\begin{equation*}
\Delta_n(y)=\prod_{1\le i<j\le n}(y_i-y_j).
\end{equation*}
\end{prop}

\begin{proof}
For a rational function $f(t)$ in variable $t$,
we write
\begin{equation*}
f(\xs)\deq\begin{pmatrix}
f(x_1) \\ \vdots \\ f(x_{kn})
\end{pmatrix}\in\Mat_{kn,1}.
\end{equation*}
Using this convention, we have
\begin{equation*}
C_{n,k}(x,y)=\kakko{\frac1{\xs+y_1},\dots,\frac1{\xs+y_n}},\quad
V_{n,k}(x)=(\xs^{n-1},\dots,\xs,1).
\end{equation*}

By Lemma \ref{lem:key_lemma_2}, we have
\begin{equation*}
\begin{split}
&\wrdet\kakko{\frac1{\xs+y_1},\dots,\frac1{\xs+y_n}}\\
=&\wrdet\kakko{\frac1{\xs+y_1},\frac1{\xs+y_2}-\frac1{\xs+y_1},\dots,\frac1{\xs+y_n}-\frac1{\xs+y_1}}\\
=&\wrdet\kakko{\frac1{\xs+y_1},\frac{y_1-y_2}{(\xs+y_1)(\xs+y_2)},\dots,\frac{y_1-y_n}{(\xs+y_1)(\xs+y_n)}}\\
=&(y_1-y_2)^k\dotsb(y_1-y_n)^k
\wrdet\kakko{\frac1{\xs+y_1},\frac1{(\xs+y_1)(\xs+y_2)},\dots,\frac1{(\xs+y_1)(\xs+y_n)}}.
\end{split}
\end{equation*}
Iterating this procedure, we reach to the expression
\begin{equation*}
\begin{split}
&\wrdet\kakko{\frac1{\xs+y_1},\dots,\frac1{\xs+y_n}}\\
=&\Delta_n(y)^k
\wrdet\kakko{\frac1{\xs+y_1},\frac1{(\xs+y_1)(\xs+y_2)},
\dots,\prod_{j=1}^n\frac1{(\xs+y_j)}}.
\end{split}
\end{equation*}
Using the multilinearity of $\kdet$ with respect to the \emph{row} vectors,
we have
\begin{equation*}
\begin{split}
&\wrdet\kakko{\frac1{\xs+y_1},\frac1{(\xs+y_1)(\xs+y_2)},
\dots,\prod_{j=1}^n\frac1{(\xs+y_j)}}\\
=&\multprod{1\le i\le kn \\ 1\le j\le n}\frac1{x_i+y_j}
\wrdet\kakko{
\prod_{j=2}^n(\xs+y_j),\prod_{j=3}^n(\xs+y_j),\dots,(\xs+y_n),1}.
\end{split}
\end{equation*}
The last wreath determinant is equal to
$\wrdet(\xs^{n-1},\dots,\xs,1)=\wrdet V_{n,k}(x)$
by Lemma \ref{lem:key_lemma_2}.
This completes the proof.
\end{proof}

We note that the proof above is exactly a wreath-analogue
of the one of the Cauchy formula \cite{Weyl}.

\begin{ex}[$k=1$]
When $k=1$, formula \eqref{eq:Cauchy-Vandermonde} is nothing but the
ordinary Cauchy determinant formula
\begin{equation*}
\det\kakko{\frac1{x_i+y_j}}_{1\le i,j\le n}
=\frac{\Delta_n(x)\Delta_n(y)}{\prod_{i,j=1}^n(x_i+y_j)}.
\end{equation*}
\end{ex}

\begin{ex}[$k=2$]
When $k=2$, \eqref{eq:Cauchy-Vandermonde} gives the formula
\begin{equation*}
\begin{split}
\wrdet[2]
\begin{pmatrix}
\dfrac1{x_1+y_1} & \dfrac1{x_1+y_2} & \dots & \dfrac1{x_1+y_n} \\
\dfrac1{x_2+y_1} & \dfrac1{x_2+y_2} & \dots & \dfrac1{x_2+y_n} \\
\vdots & \vdots & \ddots & \vdots \\
\dfrac1{x_{2n}+y_1} & \dfrac1{x_{2n}+y_2} &
\dots & \dfrac1{x_{2n}+y_n}
\end{pmatrix}
=\frac{\displaystyle\prod_{1\le i<j\le n}(y_i-y_j)^2}
{\displaystyle\multprod{1\le i\le 2n\\1\le j\le n}(x_i+y_j)}
\wrdet[2]
\begin{pmatrix}
x_1^{n-1} & \dots & x_1 & 1 \\
x_2^{n-1} & \dots & x_2 & 1 \\
\vdots & \ddots & \vdots & \vdots \\
x_{2n}^{n-1} & \dots & x_{2n} & 1
\end{pmatrix}.
\end{split}
\end{equation*}
\end{ex}

We notice that the other variant of this Cauchy-type identity also follows
immediately from \eqref{eq:Cauchy-Vandermonde}.
Indeed, we have
\begin{equation*}
\begin{split}
\wrdet\kakko{\frac1{1-\xs y_1},\dots,\frac1{1-\xs y_n}}
=\frac{\Delta_n(y)^k}
{\multprod{1\le i\le kn \\ 1\le j\le n}(1-x_iy_j)}
\wrdet V_{n,k}(x),
\end{split}
\end{equation*}
which is a wreath determinant analogue of the formula
\begin{equation*}
\det\kakko{\frac1{1-x_iy_j}}_{1\le i,j\le n}
=\frac{\Delta_n(x)\Delta_n(y)}{\prod_{i,j=1}^n(1-x_iy_j)}.
\end{equation*}

As a corollary of Theorem \ref{thm:detexp_of_wrdet},
we have the
\begin{thm}\label{thm:detexp_of_wrvdm}
The wreath Vandermonde determinant $\wrdet V_{n,k}(x)$ is given by
\begin{equation*}
\wrdet V_{n,k}(x)=
\sum_{T\in\STab((k^n))}
\wrdet I(T)\cdot\Delta_T(x),
\end{equation*}
where $\Delta_T(x)$ is the Specht polynomial
for a standard tableau $T=(t_{ij})\in\STab((k^n))$
defined by the product
\begin{equation*}
\Delta_T(x)\deq\prod_{i=1}^k\Delta_n(x_{t_{1i}},\dots,x_{t_{ni}})
\end{equation*}
of difference products.
\qed
\end{thm}

Another (symmetric) expression for $\wrdet V_{n,k}(x)$
also follows from
Corollary \ref{cor:symmetric_exp_of_wrdet}.
\begin{thm}\label{thm:symmetric_sum_of_wdet}
The equality
\begin{equation*}
\wrdet V_{n,k}(x)
=\frac1{k^{kn}}\sum_{\sigma\in S_k^n}\sigma\cdot\Delta_{n,k}(x)
\end{equation*}
holds
where $\Delta_{n,k}(x)$ is given by
\begin{equation*}
\Delta_{n,k}(x)\deq\prod_{l=1}^k \Delta_n(x_l,x_{l+k},\dots,x_{l+(n-1)k})
=\Delta_{T_0}(x).
\end{equation*}
\qed
\end{thm}

For a partition $\lambda=(\lambda_1,\dots,\lambda_N)$ of depth at most $N$,
the Schur function $s_\lambda(x_1,\dots,x_N)$ of $N$ variables
is defined as the ratio of the Vandermonde-type determinants as 
\begin{equation*}
s_\lambda(x_1,\dots,x_N)
=\dfrac{\det\kakko{x_i^{\lambda_j+N-j}}_{1\le i,j\le N}}
{\det\kakko{x_i^{N-j}}_{1\le i,j\le N}}. 
\end{equation*}
An arbitrary symmetric function can be written as a
linear combination of the Schur functions. 
We show that any symmetric function in $kn$ variables
can be written as a linear combination of the ratios
of the Vandermonde type $-\frac1k$-determinants analogously.

We recall the Cauchy identity
concerning the Schur functions (see, e.g. \cite{Mac, Weyl}).
\begin{lem}
For $m, n\in\N$, the equality
\begin{equation*}
\multprod{1\le i\le m\\ 1\le j\le n}\frac1{1-x_iy_j}\label{eq:CD-variant}
=\sum_{\len\lambda\le \min\{m,n\}}s_\lambda(x_1,\dots,x_m)s_\lambda(y_1,\dots,y_n)
\end{equation*}
holds.
\qed
\end{lem}

By the multilinearity of $\kdet$ with respect to \emph{column} vectors,
we have the following expansion formula
\begin{equation*}
\begin{split}
&\wrdet\kakko{\frac1{1-\xs y_1},\dots,\frac1{1-\xs y_n}}\\
=&\kdet{\kakko{\sum_{i_{11}\ge0}(\xs y_1)^{i_{11}},
\sum_{i_{21}\ge0}(\xs y_2)^{i_{21}},\dots,
\sum_{i_{kn}\ge0}(\xs y_n)^{i_{kn}}}}\\
=&\sum_{i_{11},i_{21},\dots,i_{kn}\ge0}
y_1^{i_{11}+\dotsb+i_{k1}}\dotsb y_n^{i_{1n}+\dotsb+i_{kn}}
\kdet\kakko{\xs^{i_{11}},\xs^{i_{21}},\dots,\xs^{i_{kn}}}.
\end{split}
\end{equation*}
Thus we have
\begin{equation}\label{eq:basic-formula}
\begin{split}
\sum_{i_{11},i_{21},\dots,i_{kn}\ge0}
y_1^{i_{11}+\dotsb+i_{k1}}\dotsb y_n^{i_{1n}+\dotsb+i_{kn}}
\kdet\kakko{\xs^{i_{11}},\xs^{i_{21}},\dots,\xs^{i_{kn}}}\\
=\Delta_n(y)^k\wrdet V_{n,k}(x)
\sum_{\len\lambda\le n}s_\lambda(x_1,\dots,x_{kn})s_\lambda(y_1,\dots,y_n).
\end{split}
\end{equation}
Comparing the homogeneous terms in \eqref{eq:basic-formula},
we have the
\begin{lem}
Put
\begin{equation*}
H_{n,k}^d(x,y)\deq
\multsum{i_{11},i_{21},\dots,i_{kn}\ge0\\ i_{11}+\dotsb+i_{kn}=d+\frac{kn(n-1)}2}
y_1^{i_{11}+\dotsb+i_{k1}}\dotsb y_n^{i_{1n}+\dotsb+i_{kn}}
\kdet\kakko{\xs^{i_{11}},\xs^{i_{21}},\dots,\xs^{i_{kn}}}.
\end{equation*}
Then, the equalities
\begin{equation*}
\wrdet\kakko{\frac1{1-\xs y_1},\dots,\frac1{1-\xs y_n}}=
\sum_{d=0}^\infty H_{n,k}^d(x,y)
\end{equation*}
and
\begin{equation*}
H_{n,k}^d(x,y)=\Delta_n(y)^k\wrdet V_{n,k}(x)
\multsum{\len\lambda\le n\\\abs\lambda=d}s_\lambda(x)s_\lambda(y)
\end{equation*}
hold.
\qed
\end{lem}

Since the Schur functions of $n$ variables are
the irreducible characters of the unitary group $U(n)$,
it follows from \eqref{eq:basic-formula} that
\begin{equation*}
\begin{split}
&s_\lambda(x_1,\dots,x_{kn})\\
=&\multsum{i_{11},i_{21},\dots,i_{kn}\ge0\\
i_{11}+\dotsb+i_{kn}=\abs\lambda+\frac{kn(n-1)}2}
\ckakko{\int_{T_n}\!
\frac{y_1^{i_{11}+\dotsb+i_{k1}}\dotsb y_n^{i_{1n}+\dotsb+i_{kn}}s_\lambda(y)}
{\Delta_n(y)^k}dg(y)}
\frac{\kdet\kakko{\xs^{i_{11}},\xs^{i_{21}},\dots,\xs^{i_{kn}}}}
{\wrdet V_{n,k}(x)},
\end{split}
\end{equation*}
where $T_n$ is the $n$-torus in $U(n)$ and $dg$ is its normalized Haar measure.
Thus implicitly, we find 
the Schur function
$s_\lambda(x_1,\dots,x_{kn})$ can be written as a
linear combination of the ratios
${\kdet\kakko{\xs^{i_{11}},\xs^{i_{21}},\dots,\xs^{i_{kn}}}}\big/
{\wrdet V_{n,k}(x)}$
of Vandermonde type $-\frac1k$-determinants.
Actually, we have the following expression.
\begin{prop}\label{thm:schur_by_k-det}
For a given sequence $\va=(a_1,\dots,a_{kn})\in\Z_{\ge0}^{kn}$ of
non-negative integers,
put
\begin{equation*}
D_{n,k}(x;\va)=\kdet\kakko{x_i^{a_j}}_{1\le i,j\le kn}.
\end{equation*}
Let us also define $\ve_i,\vdelta_{n,k}\in\Z_{\ge0}^{kn}$ by
\begin{equation*}
\ve_i=(0,\dots,0,\overset{i}{1},0,\dots,0),\quad
\vdelta_{n,k}=\sum_{j=1}^{kn}\kakko{n-1-\floor{\frac{j-1}k}}\ve_j.
\end{equation*}
Then, the Schur function $s_\lambda(x)$ is written as 
\begin{equation*}
\begin{split}
s_\lambda(x)
&=\frac1{\wrdet V_{n,k}(x)}\multsum{\mu\le\lambda \\ \abs\mu=\abs\lambda}
\sum_{\sigma\in\sym{kn}/\sym\mu}
K_{\lambda\mu}\cdot
D_{n,k}\kakko[2]{x;\vdelta_{n,k}+\sum_{i=1}^{kn}\mu_{\sigma(i)}\ve_i}.
\end{split}
\end{equation*}
\end{prop}

We notice that $\wrdet V_{n,k}(x)=D_{n,k}(x;\vdelta_{n,k})$.

For a partition $\lambda=(\lambda_1,\dots,\lambda_{kn})$ whose depth is at most $kn$,
the monomial symmetric function $m_\lambda(x)$ is defined by
\begin{equation*}
m_\lambda(x)
=\sum_{\sigma\in\sym{kn}/\sym\lambda}\prod_{i=1}^{kn}x_i^{\lambda_{\sigma(i)}}. 
\end{equation*}
Here $\sym\lambda$ is the stabilizer of $\lambda$,
that is,
$\sym\lambda=\set{\sigma\in\sym{kn}}{\lambda_{\sigma(i)}
=\lambda_i,\ 1\le i\le kn}$.
The proposition follows from the following simple lemma.
\begin{lem}
Let $\lambda$ be a partition whose depth is at most $kn$.
Then, the monomial symmetric function $m_\lambda(x)$ has the following expression
\begin{equation*}
m_\lambda(x)=\frac1{\wrdet V_{n,k}(x)}
\sum_{\sigma\in\sym{kn}/\sym\lambda}
D_{n,k}\kakko[2]{x;\vdelta_{n,k}+\sum_{i=1}^{kn}\lambda_{\sigma(i)}\ve_i}.
\end{equation*}
\end{lem}

\begin{proof}
For any $\sigma\in\sym{kn}$, we have
\begin{equation*}
\begin{split}
D_{n,k}\kakko[2]{x;\vdelta_{n,k}+\sum_{i=1}^{kn}\lambda_{\sigma(i)}\ve_i}
&=\sum_{\tau\in\sym{kn}}\kakko{-\frac1k}^{kn-\nu_{kn}(\tau)}
\prod_{i=1}^{kn}x_{\tau(i)}^{n-1-\floor{\frac{i-1}k}}\cdot
\prod_{i=1}^{kn}x_{\tau(i)}^{\lambda_{\sigma(i)}}.
\end{split}
\end{equation*}
Hence it follows that
\begin{equation*}
\begin{split}
\sum_{\sigma\in\sym{kn}}
D_{n,k}\kakko[2]{x;\vdelta_{n,k}+\sum_{i=1}^{kn}\lambda_{\sigma(i)}\ve_i}
&=\sum_{\tau\in\sym{kn}}\kakko{-\frac1k}^{kn-\nu_{kn}(\tau)}
\prod_{i=1}^{kn}x_{\tau(i)}^{n-1-\floor{\frac{i-1}k}}\cdot
\kakko{\sum_{\sigma\in\sym{kn}}\prod_{i=1}^{kn}x_{\tau(i)}^{\lambda_{\sigma(i)}}}\\
&=\wrdet V_{n,k}(x)\card{\sym\lambda}m_\lambda(x).
\end{split}
\end{equation*}
Therefore we obtain
\begin{equation*}
\begin{split}
m_\lambda(x)
&=\frac1{\wrdet V_{n,k}(x)}\frac1{\card{\sym\lambda}}
\sum_{\sigma\in\sym{kn}}
D_{n,k}\kakko[2]{x;\vdelta_{n,k}+\sum_{i=1}^{kn}\lambda_{\sigma(i)}\ve_i}\\
&=\frac1{\wrdet V_{n,k}(x)}
\sum_{\sigma\in\sym{kn}/\sym\lambda}
D_{n,k}\kakko[2]{x;\vdelta_{n,k}+\sum_{i=1}^{kn}\lambda_{\sigma(i)}\ve_i}.
\end{split}
\end{equation*}
This completes the proof.
\end{proof}

Since the Schur functions are written as a linear combination
\begin{equation*}
s_\lambda(x)
=\multsum{\mu\le\lambda \\ \abs\mu=\abs\lambda}K_{\lambda\mu}m_\mu(x)
\end{equation*}
of monomial symmetric functions, Proposition \ref{thm:schur_by_k-det}
follows immediately.

\begin{cor}
The power-sum symmetric functions $p_d(x)$,
the complete symmetric functions $h_d(x)$
and the elementary symmetric functions $e_d(x)$
are expressed as
\begin{equation*}
\begin{split}
p_d(x)&=\frac1{\wrdet V_{n,k}(x)}\sum_{i=1}^{kn}
D_{n,k}(x;\vdelta_{n,k}+d\ve_i),\\
h_d(x)&=\frac1{\wrdet V_{n,k}(x)}\sum_{1\le i_1\le\dots\le i_d\le kn}
D_{n,k}\kakko[2]{x;\vdelta_{n,k}+\sum_{j=1}^d\ve_{i_j}},\\
e_d(x)&=\frac1{\wrdet V_{n,k}(x)}\sum_{1\le i_1<\dots<i_d\le kn}
D_{n,k}\kakko[2]{x;\vdelta_{n,k}+\sum_{j=1}^d\ve_{i_j}}.
\end{split}
\end{equation*}
\qed
\end{cor}

\section{Generalities on $(n,k)$-sign and spherical functions}
\label{sec:sgn}

For $k,n\in\N$, we put
\begin{equation*}
\rnk_{n,k}\deq\set[2]{f:[kn]\to[n]}{\card{f^{-1}(j)}=k,\ \forall j\in[n]}.
\end{equation*}
We notice that $\rnk_{n,1}=\sym n$.
We also notice that $\sym{kn}$ acts on $\rnk_{n,k}$ transitively from the right,
and $\sym n$ acts on $\rnk_{n,k}$ from the left.

For $f\in\rnk_{n,k}$, we define the \emph{$(n,k)$-sign} of $f$ by
\begin{equation*}
\ksgn(f)\deq\wrdet(\delta_{f(i),j})_{\substack{1\le i\le kn\\ 1\le j\le n}}.
\end{equation*}
We see that
\begin{equation*}
\ksgn(\tau\cdot f)=\sgn(\tau)^k\ksgn(f)
\end{equation*}
for $\tau\in\sym n$.
Using this sign for $f\in\rnk_{n,k}$ and 
the very definition (4.6) of the wreath determinant we have the
\begin{lem}\label{eq:monomial_expansion_of_wrdet}
Let $k,n\in\N$. Then the equality
\begin{equation*}
\wrdet A=\sum_{f\in\rnk_{n,k}}\ksgn(f)\prod_{i\in[kn]}a_{if(i)}
\end{equation*}
holds for any $A=(a_{ij})\in\Mat_{kn,n}$.
\qed
\end{lem}

We define the element $\iota_{n,k}\in\rnk_{n,k}$ by
\begin{equation*}
\iota_{n,k}((i-1)k+j)=i\quad(1\le i\le n,\,1\le j\le k).
\end{equation*}
The stabilizer of $\iota_{n,k}$ in $\sym{kn}$ is $S_k^n$.
Hence, it follows that
\begin{equation*}
\begin{split}
\ksgn(f)
&=\sum_{w\in\sym{kn}}\kakko{-\frac1k}^{kn-\nu_{kn}(w)}
\prod_{i=1}^n\prod_{j=1}^k \delta_{fw((i-1)k+j),i}\\
&=\sum_{w\in\sym{kn}}\kakko{-\frac1k}^{kn-\nu_{kn}(w)} \delta_{fw,\iota_{n,k}}\\
&=\sum_{w\in S_k^n}\kakko{-\frac1k}^{kn-\nu_{kn}(g(f)w)},
\end{split}
\end{equation*}
where $g(f)\in\sym{kn}$ is defined by $f=\iota_{n,k}\cdot g(f)$.
Therefore, if we regard a standard tableau $T=(t_{ij})\in\STab((k^n))$
as an element of $\rnk_{n,k}$ by the assignment $T:[kn]\ni t_{ij}\mapsto i\in[n]$,
then $\ksgn(T)=\wrdet I(T)$. Hence,
the result of Theorem \ref{thm:detexp_of_wrdet} can be expressed also as
\begin{equation*}
\wrdet A=\sum_{T\in\STab((k^n))} \ksgn(T) \tdet_T(A).
\end{equation*}

We consider the injection
\begin{equation*}
\omega:\sym{n}^k\ni(w_1,\dots,w_k)\longmapsto
\kakko[1]{(i-1)k+j\mapsto w_j(i)}\in\rnk_{n,k},
\end{equation*}
and denote its image by $\snk_{n,k}$.
By Lemmas \ref{lem:double_multiplex} and \ref{eq:monomial_expansion_of_wrdet},
we have
\begin{equation}\label{eq:sgn_to_sgn}
\kakko{\frac{k!}{k^k}}^{\!\!n}
\sum_{w\in\sym n^k}\sgn(w)
\prod_{i=1}^n\prod_{j=1}^k a_{i,\omega(w)((i-1)k+j)}
=\sum_{f\in\rnk_{n,k}}\ksgn(f)
\prod_{i=1}^n\prod_{j=1}^n a_{i,f((i-1)k+j)}
\end{equation}
for $(a_{ij})_{1\le i,j\le n}\in\Mat_n$.
Comparing the coefficients in both sides,
we obtain the
\begin{cor}
For any $f\in\rnk_{n,k}$,
the equality
\begin{equation*}
\ksgn(f)=\sgn(w)\kakko{\frac{k!}{k^k}}^{\!\!n}
\frac{\card{\kakko{f\cdot S_k^n}\cap \snk_{n,k}}}{\card{f\cdot S_k^n}}
\end{equation*}
holds for $w\in\sym n^k$ such that
$\omega(w)\in \kakko{f\cdot S_k^n}\cap \snk_{n,k}$.
The sign $\sgn(w)$ does not depend on the choice of $w$.
\end{cor}

\begin{proof}
Fix an element $f\in\rnk_{n,k}$.
We notice that
the monomial $\prod_{i=1}^n\prod_{j=1}^n a_{i,f((i-1)k+j)}$
in the right-hand side of \eqref{eq:sgn_to_sgn} depends only on
the orbit $f\cdot S_k^n$.
We also notice that
the function $\ksgn$ is constant on each $S_k^n$-orbit.
Hence the coefficient of the monomial
$\prod_{i=1}^n\prod_{j=1}^n a_{i,f((i-1)k+j)}$
in the right-hand side is $\ksgn(f)\card{f\cdot S_k^n}$.
For any $w=(w_1,\dots,w_k)\in\sym{n}^k$
such that $\omega(w)\in f\cdot S_k^n$,
the sign $\sgn(w)=\sgn(w_1\dots w_k)$ gives the same value,
which can be verified by counting the inversion numbers.
It follows that
the coefficient of the monomial
$\prod_{i=1}^n\prod_{j=1}^n a_{i,f((i-1)k+j)}$
in the left-hand side is $\sgn(w)\card{\kakko{f\cdot S_k^n}\cap\snk_{n,k}}$
for any $w\in\kakko{f\cdot S_k^n}\cap\snk_{n,k}$.
Thus we have the desired conclusion.
\end{proof}

As a corollary of the discussion above,
we obtain the 
\begin{prop}
{\upshape (1)}
Put
\begin{equation*}
m_{ij}(f)=\card{\set[1]{l\in[k]}{f((i-1)k+l)=j}}.
\end{equation*}
Then
\begin{equation*}
\card{f\cdot S_k^n}=\frac{k!^n}{\prod_{i,j}m_{ij}(f)!}.
\end{equation*}
\medskip

\noindent
{\upshape (2)}
The equality
\begin{equation*}
\ksgn(f)\det(A)^k=\sum_{h\in\rnk_{n,k}}\ksgn(h)\prod_{i=1}^{kn}a_{f(i)h(i)}
\end{equation*}
holds
for any $f\in\rnk_{n,k}$ and $A=(a_{ij})_{1\le i,j\le n}\in\Mat_n$.
\pfkakko{When $k=1$, this is just the definition of the determinant.}
\medskip

\noindent
{\upshape (3)}
For $f\in\rnk_{n,k}$, put
\begin{equation*}
P_f(x_{11},\dots,x_{nk})\deq
\frac1{\card{\sym k^n}}\sum_{(\sigma_1,\dots,\sigma_n)\in\sym k^n}
\prod_{i=1}^n\prod_{j=1}^k
x_{f((i-1)k+j),\sigma_i(j)}.
\end{equation*}
Then
\begin{equation*}
\frac{\card{\kakko{f\cdot S_k^n}\cap\snk_{n,k}}}{\card{f\cdot S_k^n}}
=\text{\upshape
the coefficient of $\displaystyle\multprod{1\le i\le n\\ 1\le j\le k}x_{ij}$
in $P_f(x_{11},\dots,x_{nk})$}.
\end{equation*}
\qed
\end{prop}

It is convenient to express an element $f\in\rnk_{n,k}$
as an $n\times k$ matrix
whose $(i,j)$-entry is given by $f((i-1)k+j)$, that is,
\begin{equation*}
f=\begin{pmatrix}
f(1) & \dots & f(k) \\
\vdots & \ddots & \vdots \\
f((n-1)k+1) & \dots & f(nk)
\end{pmatrix}.
\end{equation*}
If $f_1,f_2\in\rnk_{n,k}$ and
$f_2=f_1\cdot\sigma$ for some $\sigma\in S_k^n$,
then each row vector of ${f_2}$ is a permutation
of the corresponding row vector of ${f_1}$.

\begin{ex}
Let us calculate $\ksgn(U_4)=\wrdet[2]I(U_4)$ for $U_4$
(regarding as an element in $\rnk_{3,2}$)
given in Example \ref{ex:detexp_of_wrdet_for_(n,k)=(3,2)}.
In the matrix notation,
\begin{equation*}
U_4=\young(13,25,46)=
\ckakko{
\begin{matrix}
1 \mapsto 1 & 2 \mapsto 2 \\
3 \mapsto 1 & 4 \mapsto 3 \\
5 \mapsto 2 & 6 \mapsto 3
\end{matrix}
}
=\begin{pmatrix}
1 & 2 \\ 1 & 3 \\ 2 & 3
\end{pmatrix}.
\end{equation*}
It follows that
\begin{equation*}
U_4\cdot S_2^3=
\ckakko{
\begin{pmatrix}
1 & 2 \\ 1 & 3 \\ 2 & 3
\end{pmatrix},\,
\begin{pmatrix}
1 & 2 \\ 1 & 3 \\ 3 & 2
\end{pmatrix},\,
\begin{pmatrix}
1 & 2 \\ 3 & 1 \\ 2 & 3
\end{pmatrix},\,
\begin{pmatrix}
1 & 2 \\ 3 & 1 \\ 3 & 2
\end{pmatrix},\,
\begin{pmatrix}
2 & 1 \\ 1 & 3 \\ 2 & 3
\end{pmatrix},\,
\begin{pmatrix}
2 & 1 \\ 1 & 3 \\ 3 & 2
\end{pmatrix},\,
\begin{pmatrix}
2 & 1 \\ 3 & 1 \\ 2 & 3
\end{pmatrix},\,
\begin{pmatrix}
2 & 1 \\ 3 & 1 \\ 3 & 2
\end{pmatrix}
}
\end{equation*}
and
\begin{equation*}
\kakko{
U_4\cdot S_2^3}\cap\snk_{3,2}=
\ckakko{
\begin{pmatrix}
1 & 2 \\ 3 & 1 \\ 2 & 3
\end{pmatrix},\,
\begin{pmatrix}
2 & 1 \\ 1 & 3 \\ 3 & 2
\end{pmatrix}
}.
\end{equation*}
Since
\begin{equation*}
\begin{pmatrix}
1 & 2 \\ 3 & 1 \\ 2 & 3
\end{pmatrix}=\omega((2,3),(1,2))
\end{equation*}
and $\sgn((2,3),(1,2))=1$
(where $(i,j)$ denotes the transposition of $i$ and $j$), we get
\begin{equation*}
\wrdet[2]I(U_4)=\kakko{\frac{2!}{2^2}}^{\!\!3}\times\frac28=\frac1{32}.
\end{equation*}
We remark that
\begin{equation*}
\begin{pmatrix}
m_{11}(U_4) & m_{12}(U_4) & m_{13}(U_4) \\
m_{21}(U_4) & m_{22}(U_4) & m_{23}(U_4) \\
m_{31}(U_4) & m_{32}(U_4) & m_{33}(U_4)
\end{pmatrix}
=\begin{pmatrix}
1 & 1 & 0 \\ 1 & 0 & 1 \\ 0 & 1 & 1
\end{pmatrix}
\end{equation*}
and we see that
\begin{equation*}
\frac{2!^3}{\prod_{1\le i,j\le 3}m_{ij}(U_4)!}
=\frac{2!^3}{1!1!0!1!0!1!0!1!1!}=8=\card{U_4\cdot S_2^3}
\end{equation*}
as we counted above.
We also note that
\begin{equation*}
\begin{split}
P_{U_4}(x_{11},\dots,x_{32})
&=\frac18
(x_{11}x_{22}+x_{12}x_{21})(x_{11}x_{32}+x_{12}x_{31})(x_{21}x_{32}+x_{22}x_{31}),
\end{split}
\end{equation*}
and the coefficient of $x_{11}x_{21}x_{31}x_{12}x_{22}x_{32}$
of $P_{U_4}$ is $\dfrac28=\dfrac14$.
\end{ex}

Let us put
\begin{equation}\label{eq:def_of_varphi}
\varphi_{n,k}(g)=\frac{\kdet\kakko{g\cdot\1_k^{\oplus n}}}
{\kdet\kakko{\1_k^{\oplus n}}}
=k^{kn}\frac1{\card{S_k^n}}
\sum_{\sigma\in S_k^n}\kakko{-\frac1k}^{kn-\nu_{kn}(g^{-1}\sigma)}
\end{equation}
for $g\in\sym{kn}$.
We note that $\varphi_{n,k}(g^{-1})=\varphi_{n,k}(g)$
since $\nu_{kn}(g^{-1}\sigma)=\nu_{kn}(g\sigma^{-1})$.
By Lemma \ref{lem:relative_invariance_of_wrdet}
and its $S_k^n$-invariance of $\1_k^{\oplus n}$,
it follows that
\begin{equation*}
\varphi_{n,k}(h_1gh_2)=\chi_{n,k}(h_1h_2)^k\varphi_{n,k}(g)
\end{equation*}
for $g\in\sym{kn}$ and $h_1, h_2\in\wsym kn$.
In particular, $\varphi_{n,k}$ is a $S_k^n$-biinvariant
(or $S_k^n$-zonal spherical) function on $\sym{kn}$.
We note that the rightmost side of \eqref{eq:def_of_varphi}
can be considered as an analogue of
the integral expression of the zonal spherical function
of a Riemannian symmetric space due to Harish-Chandra 
(see, e.g. \cite{Helgason}).

\begin{lem}\label{lem:positivity_of_phi}
The $\chi_{n,k}^k$-spherical function $\varphi_{n,k}$ relative to  
the wreath product $\wsym kn$ on $\sym{kn}$
is expressed as a matrix element of the \pfkakko{unitary} representation 
$M_{n,k}^{T_{kn},\det} (\cong\Smod{kn}\lambda)$ of $\sym{kn}$:
\begin{equation*}
\varphi_{n,k}(g)
=\frac{\inprod{g\cdot\wrdet(X)}{\wrdet(X)}}{\inprod{\wrdet(X)}{\wrdet(X)}},
\end{equation*}
where $\inprod{\,}{\,}$ denotes the invariant inner product
on $M_{n,k}^{T_{kn},\det}$.
In particular,
$\varphi_{n,k}$ is a positive definite function.
\end{lem}

\begin{proof}
Consider the projection
\begin{equation*}
P_{n,k}=\frac1{\card{S_k^n}}\sum_{\sigma\in S_k^n}\sigma\in\C[\sym{kn}].
\end{equation*}
By \eqref{eq:S_k^n-invariant_space_is_one_dimensional},
for each $g\in\sym{kn}$,
there exists a constant $C(g)$ such that
\begin{equation}\label{eq:project_into_wrdet}
P_{n,k}\,g\cdot\wrdet(X)=C(g)\wrdet(X).
\end{equation}
Since $P_{n,k}$ is self-adjoint with respect to $\inprod{\,}{}$
and $\wrdet(X)$ is $S_k^n$-invariant,
it follows that
\begin{equation*}
\begin{split}
\inprod{g\cdot\wrdet(X)}{\wrdet(X)}
&=\inprod{g\cdot\wrdet(X)}{P_{n,k}\cdot\wrdet(X)}\\
&=\inprod{P_{n,k}g\cdot\wrdet(X)}{\wrdet(X)}
=C(g)\inprod{\wrdet(X)}{\wrdet(X)}.
\end{split}
\end{equation*}
To determine $C(g)$,
let us calculate the coefficient of
$\prod_{p=1}^n\prod_{l=1}^k x_{(p-1)k+l,p}$
in the both sides of \eqref{eq:project_into_wrdet}.
It is immediate to see that the coefficient in the right-hand side
is $C(g)\kakko{\frac{k!}{k^k}}^n=C(g)\kdet\kakko{\1_k^{\oplus n}}$.
We look at the left-hand side:
\begin{equation*}
\begin{split}
P_{n,k}\,g\cdot\wrdet(X)
&=\frac1{\card{S_k^n}}\sum_{\sigma\in S_k^n}\sigma g\cdot\wrdet(X)\\
&=\frac1{\card{S_k^n}}\sum_{\sigma\in S_k^n}
\sum_{h\in\sym{kn}}\kakko{-\frac1k}^{\!kn-\nu_{kn}(h)}
\prod_{p=1}^n\prod_{l=1}^k x_{(\sigma g h)((p-1)k+l),p}\\
&=\sum_{h\in\sym{kn}}\kakko{
\frac1{\card{S_k^n}}\sum_{\sigma\in S_k^n}
\kakko{-\frac1k}^{\!kn-\nu_{kn}(g^{-1}\sigma^{-1}h)}}
\prod_{p=1}^n\prod_{l=1}^k x_{h((p-1)k+l),p}.
\end{split}
\end{equation*}
Hence the coefficient of $\prod_{p=1}^n\prod_{l=1}^k x_{(p-1)k+l,p}$
in $P_{n,k}\,g\cdot\wrdet(X)$ is equal to 
\begin{equation*}
\sum_{h\in S_k^n}\frac1{\card{S_k^n}}\sum_{\sigma\in S_k^n}
\kakko{-\frac1k}^{\!kn-\nu_{kn}(g^{-1}\sigma^{-1}h)}
=\sum_{\sigma\in S_k^n}
\kakko{-\frac1k}^{\!kn-\nu_{kn}(g^{-1}\sigma)}
=\kdet\kakko{g\cdot\1_k^{\oplus n}}.
\end{equation*}
Thus we have
\begin{equation*}
C(g)=\frac{\kdet\kakko{g\cdot\1_k^{\oplus n}}}{\kdet\kakko{\1_k^{\oplus n}}}
=\varphi_{n,k}(g).
\end{equation*}
This completes the proof.
\end{proof}

\begin{rem}
By specializing
the Frobenius character formula for $\sym N$,
we have
\begin{equation*}
\alpha^{N-\nu_N(g)}
=\sum_{\lambda\vdash N}\frac{f^\lambda}{N!}f_\lambda(\alpha)\chi^\lambda(g)
\qquad(g\in\sym N), 
\end{equation*}
where $f_\lambda(\alpha)$ denotes the content polynomial defined by
\begin{equation*}
f_\lambda(\alpha)=\prod_{(i,j)\in\lambda}(1+(j-i)\alpha).
\end{equation*}
Since
\begin{equation*}
f_\lambda\kakko{-\frac1k}
=\prod_{(i,j)\in\lambda}\kakko{1-\frac1k(j-i)}
=\frac1{k^{kn}}\prod_{(j,i)\in\lambda'}\kakko{k+(i-j)}
=\frac{(kn)!}{f^\lambda}\frac{\card{\SSTab_k(\lambda')}}{k^{kn}},
\end{equation*}
it follows that
\begin{equation*}
\kakko{-\frac1k}^{\!kn-\nu_{kn}(g)}=\sum_{\lambda\vdash kn}
\frac{\card{\SSTab_k(\lambda')}}{k^{kn}}\chi^\lambda(g).
\end{equation*}
Hence the function $\varphi_{n,k}$ is a linear combination
\begin{equation*}
\varphi_{n,k}(g)
=\sum_{\lambda\vdash kn}\card{\SSTab_k(\lambda')}
\phi_{n,k}^\lambda(g)
\end{equation*}
of $S_k^n$-zonal spherical functions
\begin{equation*}
\phi_{n,k}^\lambda(g)=\frac1{\card{S_k^n}}\sum_{\sigma\in S_k^n}
\chi^\lambda(g^{-1}\sigma)
\end{equation*}
with nonnegative (integral) coefficients.
Therefore, it is immediate to see again that
$\varphi_{n,k}$ is a positive definite function.
\end{rem}

\begin{rem}
Since
$\inprod{\ind_{S_k^n}^{\sym{kn}}\trivrpn{S_k^n}}{\Smod{kn}\lambda}
=K_{\lambda,(k^n)}$
for $\lambda\vdash kn$,
the pair $(\sym{kn},S_k^n)$
is not a Gelfand pair in general. Further, although 
one can verify that the pair $(\sym{kn},\wsym kn)$
is a Gelfand pair when $k=2$ (see p.401 in \cite{Mac}, in fact, the 
wreath product $\wsym 2n$ is isomorphic to the 
\emph{hyperoctahedral} group of degree $n$),
it is not the case for a general $k$. 
Actually, when $n=3$,
by looking at the Schur function expansion of
the plethysm $h_3\circ h_k$ (see p.141 in \cite{Mac}),
it follows that the induced representation
$\ind_{\wsym k3}^{\sym{3k}}\trivrpn{\wsym k3}$
is not multiplicity free
when $k\ge18$.
\end{rem}

For a standard tableau $T\in\STab((k^n))$, 
we define
\begin{equation*}
D_T(X)=\wrdet(g(T)^{-1}\cdot X), 
\end{equation*}
where $g(T)$ is a permutation given in \eqref{eq:def_of_g(T)}.
We see that
\begin{equation*}
\begin{split}
D_T(X)
&=\sum_{S\in\STab((k^n))}\wrdet(g(T)^{-1}I(S))\tdet_S(X) \\
&=\kakko{\frac{k!}{k^k}}^{\!n}
\sum_{S\in\STab((k^n))}\varphi_{n,k}(g(T)^{-1}g(S))\tdet_S(X).
\end{split}
\end{equation*}
We now define the $f^{(k^n)}\times f^{(k^n)}$ matrix $\Xi_{n,k}$ by 
\begin{equation}
\Xi_{n,k}=\kakko{\varphi_{n,k}(g(T)^{-1}g(S))}_{S,T\in\STab((k^n))}.
\end{equation}
Since $\varphi_{n,k}(g)=\varphi_{n,k}(g^{-1})$, one finds that 
the matrix $\Xi_{n,k}$ is symmetric. Moreover, we notice 
that $\det\Xi_{n,k}\ge0$ by Lemma \ref{lem:positivity_of_phi},
because $\varphi_{n,k}$ is a positive definite function. Then the following 
conjecture looks quite reasonable. 
\begin{conj}\label{positiveConjecture}
The matrix $\Xi_{n,k}$ is positive definite; in particular, 
one has  $\det\Xi_{n,k}>0$.
In other words, 
$\{D_T(X)\}_{T\in\STab((k^n))}$ gives
another basis of the space
$M_{n,k}^{T_{kn},\det}=\C[\sym{kn}]\cdot\wrdet$.
\end{conj}
We try to examine the first few examples which may 
support the above conjecture. 
\begin{ex}
We have
\begin{align*}
\det\Xi_{2,2}&=\frac13\kakko{\frac32}^2,\quad
\det\Xi_{3,2}=\frac23\kakko{\frac34}^5,\quad
\det\Xi_{2,3}=\frac32\kakko{\frac23}^5,\\
\det\Xi_{4,2}&=\frac{2^65}3\kakko{\frac38}^{14},\quad
\det\Xi_{2,4}=\frac3{2^65}\kakko{\frac56}^{14}.
\end{align*}
We notice here that
\begin{equation*}
f^{(2^2)}=2,\quad
f^{(2^3)}=f^{(3^2)}=5,\quad
f^{(2^4)}=f^{(4^2)}=14.
\end{equation*}
\end{ex}

\appendix

\section{Appendix : Laplace expansion of $\alpha$-determinants}
\label{sec:laplace}

\begin{prop}[Laplace expansion]
For a given $n$ by $n$ matrix $X=(x_{ij})_{1\le i,j\le n}$,
we have
\begin{equation*}
\adet X=\sum_{p=1}^n \alpha^{1-\delta_{pq}}x_{pq}\adet X_{pq},
\end{equation*}
where $X_{pq}$ is a $n-1$ by $n-1$ matrix obtained by
the following procedure:
{\upshape(1)} remove $q$-th column vector and $q$-th row vector in $X$,
{\upshape(2)} if $p\neq q$, then replace the row vector
$(x_{p1},\dots,x_{pn})$ in $X$ by $(x_{q1},\dots,x_{qn})$.
\end{prop}

\begin{proof}
We have
\begin{equation*}
\begin{split}
\adet X&=
\sum_{p=1}^n \multsum{g\in\sym n\\ g(q)=p}
\alpha^{n-\nu_n(g)}\prod_{i=1}^n x_{g(i)i}\\
&=\sum_{p=1}^n x_{pq}\multsum{g\in\sym n\\ g(q)=q}
\alpha^{n-\nu_n((p,q)\cdot g)}
\prod_{1\le i(\ne q)\le n} x_{(p,q)\cdot g(i)i}\\
&=\sum_{p=1}^n \alpha^{1-\delta_{pq}}x_{pq}
\multsum{g\in\sym n\\ g(q)=q}
\alpha^{(n-1)-\nu_{n-1}(g)}
\prod_{1\le i(\ne q)\le n} x_{(p,q)\cdot g(i)i}\\
&=\sum_{p=1}^n \alpha^{1-\delta_{pq}}x_{pq}\adet X_{pq}.
\end{split}
\end{equation*}
Here we use the fact that
$\nu_n((p,q)\cdot g)=\nu_{n-1}(g)+\delta_{pq}$
if $g(q)=q$ (see the proof of Lemma \ref{lem:shifted_cycle_sum}).
\end{proof}

\begin{ex}[$n=4$]
For $X=\begin{pmatrix}
x_{11} & x_{12} & x_{13} & x_{14} \\
x_{21} & x_{22} & x_{23} & x_{24} \\
x_{31} & x_{32} & x_{33} & x_{34} \\
x_{41} & x_{42} & x_{43} & x_{44}
\end{pmatrix}$, we have
\begin{equation*}
\begin{split}
X_{12}=
\begin{pmatrix}
x_{21} & x_{23} & x_{24} \\
x_{31} & x_{33} & x_{34} \\
x_{41} & x_{43} & x_{44}
\end{pmatrix},\quad
X_{22}=
\begin{pmatrix}
x_{11} & x_{13} & x_{14} \\
x_{31} & x_{33} & x_{34} \\
x_{41} & x_{43} & x_{44}
\end{pmatrix},\\
X_{32}=
\begin{pmatrix}
x_{11} & x_{13} & x_{14} \\
x_{21} & x_{23} & x_{24} \\
x_{41} & x_{43} & x_{44}
\end{pmatrix},\quad
X_{42}=
\begin{pmatrix}
x_{11} & x_{13} & x_{14} \\
x_{31} & x_{33} & x_{34} \\
x_{21} & x_{23} & x_{24}
\end{pmatrix}.
\end{split}
\end{equation*}
Hence we have
\begin{equation*}
\begin{split}
\adet\begin{pmatrix}
x_{11} & x_{12} & x_{13} & x_{14} \\
x_{21} & x_{22} & x_{23} & x_{24} \\
x_{31} & x_{32} & x_{33} & x_{34} \\
x_{41} & x_{42} & x_{43} & x_{44}
\end{pmatrix}
=&\alpha x_{12}
\adet\begin{pmatrix}
x_{21} & x_{23} & x_{24} \\
x_{31} & x_{33} & x_{34} \\
x_{41} & x_{43} & x_{44}
\end{pmatrix}
+x_{22}\adet
\begin{pmatrix}
x_{11} & x_{13} & x_{14} \\
x_{31} & x_{33} & x_{34} \\
x_{41} & x_{43} & x_{44}
\end{pmatrix}\\
&+\alpha x_{32}\adet
\begin{pmatrix}
x_{11} & x_{13} & x_{14} \\
x_{21} & x_{23} & x_{24} \\
x_{41} & x_{43} & x_{44}
\end{pmatrix}
+\alpha x_{42}\adet
\begin{pmatrix}
x_{11} & x_{13} & x_{14} \\
x_{31} & x_{33} & x_{34} \\
x_{21} & x_{23} & x_{24}
\end{pmatrix}.
\end{split}
\end{equation*}
\end{ex}

\smallskip

\begin{flushleft}
Kazufumi KIMOTO\\
Department of Mathematical Science,
University of the Ryukyus.\\
Senbaru, Nishihara, Okinawa 903-0231, JAPAN.\\
\texttt{kimoto@math.u-ryukyu.ac.jp}\\
\end{flushleft}

\begin{flushleft}
Masato WAKAYAMA\\
Faculty of Mathematics,
Kyushu University.\\
Hakozaki, Fukuoka 812-8518, JAPAN.\\
\texttt{wakayama@math.kyushu-u.ac.jp}
\end{flushleft}

\end{document}